\documentclass[3p,times,review]{elsarticle}

\usepackage{amssymb}

\usepackage[figuresright]{rotating}

\def\mb#1{\mbox {\boldmath {$#1$}}} 

\usepackage{stmaryrd,bm,amsmath}
\usepackage{dsfont}

\usepackage{graphicx,stmaryrd,sectsty}
\usepackage{setspace}
\usepackage[font=small,labelfont=bf]{caption}
\usepackage{epstopdf}
\usepackage{verbatim}
\usepackage{graphicx}
\usepackage{caption}
\usepackage{subcaption}
\usepackage{float}
\usepackage[toc,page]{appendix}

\newcommand{\B}{{Ba\v zant}}
\newcommand{\bc}{\begin{center}}
\newcommand{\ec}{\end{center}}
\newcommand{\bfr}{\begin{flushright}}
\newcommand{\efr}{\end{flushright}}

\newcommand{\be}{\begin{enumerate}}
\newcommand{\ee}{\end{enumerate}}
\newcommand{\bi}{\begin{itemize}}
\newcommand{\ei}{\end{itemize}}
\newcommand{\bd}{\begin{description}}
\newcommand{\ed}{\end{description}}
\newcommand{\beq}{\begin{equation}}
\newcommand{\eeq}{\end{equation}}
\newcommand{\bea}{\begin{eqnarray}}
\newcommand{\eea}{\end{eqnarray}}

\newcommand{\bfi}{\begin{figure}}
\newcommand{\efi}{\end{figure}}
\newcommand{\bay}{\begin{array}{l}}
\newcommand{\eay}{\end{array}}

\def\mb#1{\mbox {\boldmath {$#1$}}} 

\newcommand{\mbf}{\mathbf}
\newcommand*\Bell{\ensuremath{\boldsymbol\ell}}

\usepackage{tikz}
\newcommand*\circled[1]{\tikz[baseline=(char.base)]{
            \node[shape=circle,draw,inner sep=0.1pt] (char) {#1};}}

\begin{document}

\begin{titlepage}
\clearpage\thispagestyle{empty}
\noindent
\hrulefill
\begin{figure}[h!]
\centering
\includegraphics[width=2 in]{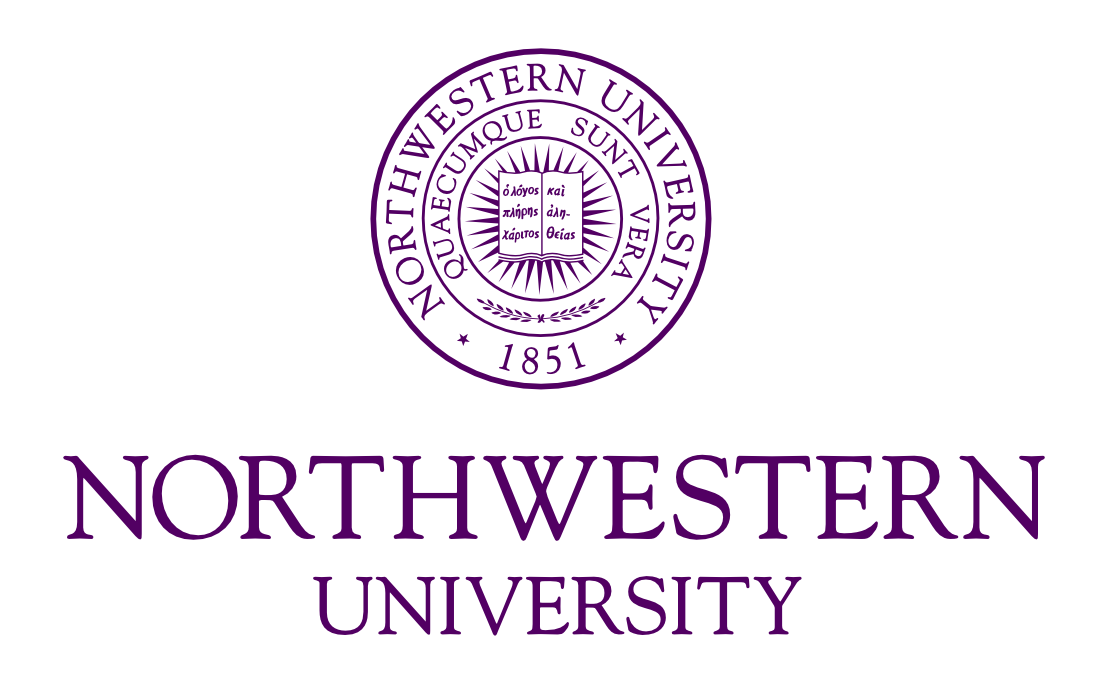}
\end{figure}
\begin{center}
{
{\bf Center for Sustainable Engineering of Geological and
Infrastructure Materials (SEGIM)} \\ [0.1in]
Department of Civil and Environmental Engineering \\ [0.1in]
McCormick School of Engineering and Applied Science \\ [0.1in]
Evanston, Illinois 60208, USA
}
\end{center} 
\hrulefill \\ \vskip 2mm
\vskip 0.5in
\begin{center}
{\large {\bf \uppercase{Discontinuous Cell Method (DCM) for the Simulation of Cohesive Fracture and Fragmentation of Continuous Media}}}\\[0.5in]
{\large {\sc Gianluca Cusatis, Roozbeh Rezakhani, Edward A. Schauffert}}\\[0.75in]
{\sf \bf SEGIM INTERNAL REPORT No. 16-08/587D}\\[0.75in]
\end{center}
\noindent {\footnotesize {{\em Submitted to Journal of Engineering Fracture Mechanics \hfill August 2016} }}
\end{titlepage}

\begin{frontmatter}




\title{Discontinuous Cell Method (DCM) for the Simulation of Cohesive Fracture and Fragmentation of Continuous Media}


\author[l1]{Gianluca Cusatis}
\author[l1]{Roozbeh Rezakhani \footnote{Corresponding author. \\ E-mail address: rrezakhani@u.northwestern.edu}}
\author[l2]{Edward A. Schauffert}
\address[l1]{Department of Civil and Environmental Engineering, Northwestern University, Evanston (IL) 60208, USA.}
\address[l2]{Greenman-Pedersen, INC., Albany (NY), 12205, USA.}

\begin{abstract}
In this paper, the Discontinuous Cell Method (DCM) is formulated with the objective of simulating cohesive fracture propagation and fragmentation in homogeneous solids without issues relevant to excessive mesh deformation typical of available Finite Element formulations. DCM discretizes solids by using the Delaunay triangulation and its associated Voronoi tessellation giving rise to a system of discrete cells interacting through shared facets. For each Voronoi cell, the displacement field is approximated on the basis of rigid body kinematics, which is used to compute a strain vector at the centroid of the Voronoi facets. Such strain vector is demonstrated to be the projection of the strain tensor at that location. At the same point stress tractions are computed through vectorial constitutive equations derived on the basis of classical continuum tensorial theories. Results of analysis of a cantilever beam are used to perform convergence studies and comparison with classical finite element formulations in the elastic regime. Furthermore, cohesive fracture and fragmentation of homogeneous solids are studied under quasi-static and dynamic loading conditions. The mesh dependency problem, typically encountered upon adopting softening constitutive equations, is tackled through the crack band approach. This study demonstrates the capabilities of DCM by solving multiple benchmark problems relevant to cohesive crack propagation. The simulations show that DCM can handle effectively a wide range of problems from the simulation of a single propagating fracture to crack branching and fragmentation.
\end{abstract}

\begin{keyword}
cohesive fracture \sep finite elements \sep discrete models \sep delaunay triangulation \sep voronoi tessellation \sep fragmentation


\end{keyword}

\end{frontmatter}

\section{Introduction}
A quantitative investigation of cohesive fracture propagation necessitates an accurate description of various fracture phenomena including: crack initiation; propagation along complex three-dimensional paths; interaction and coalescence of distributed multi-cracks into localized continuous cracks; and interaction of fractured/unfractured material.
The classical Finite Element (FE) method, although it has been used with some success to address some of these aspects, is inherently incapable of modeling the displacement discontinuities associated with fracture. To address this issue, advanced computational technologies have been developed in the recent past.
First, the embedded discontinuity methods (EDMs) were proposed to handle displacement discontinuity within finite elements. In these methods the crack is represented by a narrow band of high strain, which is embedded in the element and can be arbitrarily aligned. Many different EDM formulations can be found in the literature and a comprehensive comparative study of these formulations appears in Reference \cite{Jirasek-1}. The most common drawbacks of EDM formulations are stress locking (spurious stress transfer between the crack surfaces), inconsistency between the stress acting on the crack surface and the stress in the adjacent material bulk, and mesh sensitivity (crack path depending upon mesh alignment and refinement).

A method that does not experience stress locking and reduces mesh sensitivity is the extended finite element method (X-FEM). X-FEM, first introduced by Belytschko \& Black \cite{Belytschko-2}, exploits the partition of unity property of FE shape functions. This property enables discontinuous terms to be incorporated locally in the displacement field without the need of topology changes in the initial uncracked mesh. Mo{\"e}s et al.  \cite{Moes-1} enhanced the primary work of Belytschko et al. \cite{Belytschko-2} through including a discontinuous enrichment function to represent displacement jump across the crack faces away from the crack tip. X-FEM has been successfully applied to a wide variety of problems. Dolbow et. al. \cite{Dolbow-1} applied XFEM to the simulation of growing discontinuity in Mindlin-Reissner plates by employing appropriate asymptotic crack-tip enrichment functions. Belytschko and coworkers \cite{Belytschko-3} modeled evolution of arbitrary discontinuities in classical finite elements, in which discontinuity branching and intersection modeling are handled by the virtue of adding proper terms to the related finite element displacement shape functions. Furthermore, they studied crack initiation and propagation under dynamic loading condition and used a criterion based on the loss of hyperbolicity of the underlying continuum problem \cite{Belytschko-4}. Zi et Al. \cite{Zi-Belytschko} extended X-FEM to the simulation of cohesive crack propagation. 
The main drawbacks of X-FEM are that the implementation into existing FE codes is not straightforward, the insertion of additional degrees of freedoms is required on-the-fly to describe the discontinuous enrichment, and complex quadrature routines are necessary to integrate discontinuous integrands.

Another approach widely used for the simulation of cohesive fracture is based on the adoption of cohesive zero-thickness finite elements located at the interface between the usual finite elements that discretize the body of interest \cite{Camacho-1, Ortiz-1}. This method, even if its implementation is relatively simple, tends to be computationally intensive because of the large number of nodes that are needed to allow fracturing at each element interface. Furthermore, in the elastic phase the zero-thickness finite elements require the definition of an artificial penalty stiffness to ensure inter-element compatibility. This stiffness usually deteriorates the accuracy and rate of convergence of the numerical solution and it may cause numerical instability. To avoid this problem, algorithms have been proposed in the literature \cite{Pandolfi-1} for the dynamic insertion of cohesive fractures into FE meshes. The dynamic insertion works reasonably well in high speed dynamic applications but is not adequate for quasi-static applications and leads to inaccurate stress calculations along the crack path.

An attractive alternative to the aforementioned approaches is the adoption of discrete models (particle and lattice models), which replace the continuum a priori by a system of rigid particles that interact by means of linear/nonlinear springs or by a grid of beam-type elements. These models were first developed to describe the behavior of particulate materials \cite{Cundall-1} and to solve elastic problems in the pre-computers era \cite{Hrennikoff-1}. Later, they have been adapted to simulate fracture and failure of quasi-brittle materials in both two \cite{Schlangen-1} and three dimensional problems \cite{Cusatis-1, Cusatis-2, Lilliu-1, Cusatis-3}. In this class of models, it is worth mentioning the rigid-body-spring model developed by Bolander and collaborators, which dicretizes the material domain using Voronoi diagrams with random geometry, interconnected by zero-size springs, to simulate cohesive fracture in two and three dimensional problems \cite{Bolander-2, Bolander-1, Bolander-3, Bolander-4}. Various other discrete models, in the form of either lattice or particle models, have been quite successful recently in simulating concrete materials \cite{cusatis-ldpm-1, cusatis-ldpm-2,Rezakhani-JMPS,Leite-1,Donze-1,Grassl-1}. 

Discrete models can realistically simulate fracture propagation and fragmentation without suffering from the aforementioned typical drawbacks of other computational technologies. The effectiveness and the robustness of the method are ensured by the fact that: a) their kinematics naturally handle displacement discontinuities; b) the crack opening at a certain point depends upon the displacements of only two nodes of the mesh; c) the constitutive law for the fracturing behavior is vectorial; d) remeshing of the material domain or inclusion of additional degrees of freedom during the fracture propagation process is not necessary. Despite these advantages the general adoption of these methods to simulate fracture propagation in continuous media has been quite limited because of various drawbacks in the uncracked phase, including: 1) the stiffness of the springs is defined through a heuristic (trial-and-error) characterization; 2) various elastic phenomena, e.g. Poisson's effect, cannot be reproduced exactly; 3) the convergence of the numerical scheme to the continuum solution cannot be proved; 4) amalgamation with classical tensorial constitutive laws is not possible; and 5) spurious numerical heterogeneity of the response (not related to the internal structure of the material) is inherently associated with these methods if simply used as discretization techniques for continuum problems.

The Discontinuous Cell Method (DCM) presented in this paper provides a framework unifying discrete models and continuum based methods. The Delaunay triangulation is employed to discretize the solid domain into triangular elements, the Voronoi tessellation is then used to build a set of discrete polyhedral cells whose kinematics is described through rigid body motion typical of discrete models. Tonti  \cite{tonti-1} presented a somewhat similar approach to discretize the material domain and to compute the finite element nodal forces using dual cell geometries. Furthermore, the DCM formulation is similar to that of the discontinuous Galerkin method which has primarily been applied in the past to the solution of fluid dynamics problems, but has also been extended to the study of elasticity \cite{Guzey-1}. Recently, discontinuous Galerkin approaches have also been used for the study of fracture mechanics \cite{Shen-1} and cohesive fracture propagation \cite{Abedi-1}. The DCM formulation can be considered as a discontinuous Galerkin approach which utilizes piecewise constant shape functions. Another interesting feature of DCM is that the formulation includes rotational degrees of freedom. Researchers have attempted to introduce rotational degrees of freedom to classical finite elements by considering special form of displacement functions along each element edge to improve their performance in bending problems \cite{Allman-1,Bergan-1}. This strategy leads often to zero energy deformation modes and to singular element stiffness matrix even if the rigid body motions are constrained. DCM formulation simply incorporates nodal rotational degrees of freedom, without suffering from the aforementioned problem.

\section{Governing Equations} \label{2}
Equilibrium, compatibility, and constitutive laws for Cauchy continua can be formulated as limit case of the governing equations for Cosserat continua in which displacements and rotations are assumed to be independent fields but the couple stress tensor is identically zero \cite{Zhou-1}. For small deformations and for any position vector $\bf x$ in the material domain $\Omega$, one has 
\begin{equation}
\gamma_{ij} =u_{j,i} -e_{ijk} \varphi_k \label{Eq:Compatibility-1} 
\end{equation}
and 
\begin{equation}
\sigma_{ji,j} + b_{0i} = \rho \ddot{u}_i ;~~~ e_{ijk} \sigma_{jk} = 0 \label{Eq:Equilibrium-2}
\end{equation}

In the above equations, the summation rule of repeated indices applies; $u_i$ and $\varphi_i$ are displacement and rotation fields, respectively. $\gamma_{ij}$ is the strain tensor; $\sigma_{ij}$ is the stress tensor; $b_{0i}$ are body forces per unit volume; $\rho$ is the mass density. Subscripts $i$, $j$, and $k$ represent components of Cartesian coordinate system which can be $i,j,k = 1,2,3$ in three dimensional problems; $e_{ijk}$ is the Levi-Civita permutation symbol. Considering any position dependent field variable such as $f(\bf x, t)$, $f_{,i}$ represent partial derivative of $f$ with respect to the $i$th component of the coordinate system, while $\dot{f}$ is the time derivative of the variable. The partial differential equations above need to be complemented by appropriate boundary conditions that can either involve displacements, $u_i - u_{0i} = 0$ on $\Gamma_u$ (essential boundary conditions); or tractions, $\sigma_{ji}n_j - t_{0i}=0$ on  $\Gamma_t$ (natural boundary conditions); where $\Gamma = \Gamma_t \cup \Gamma_u$ is the boundary of the solid volume $\Omega$.

In the elastic regime, the constitutive equations can be written as
\begin{gather}\label{eq:const}
\sigma_{ij} = E_V \epsilon_V \delta_{ij} + E_D (\gamma_{ij}-\epsilon_V \delta_{ij})
\end{gather}
where $\epsilon_V=\gamma_{ii}/3$ is the volumetric strain; $E_V$ and $E_D$ are the volumetric and deviatoric moduli that can be expressed through Young's modulus $E$ and Poisson's ratio $\nu$: $E_V = E/(1-2\nu)$; $E_D=E/(1+\nu)$. It is worth observing that since the solution of the problem formulated above requires the stress tensor to be symmetric (see second equation in Equation \ref{Eq:Equilibrium-2}), the constitutive equations imply the symmetry of the strain tensor as well, which, in turn leads to displacements and rotations to be related through the following expression: $e_{ijk} \varphi_k = (u_{j,i}-u_{i,j})/2$

The weak form of the equilibrium equations can be obtained through the Principle of Virtual Work (PVW) as
\begin{equation} \label{PVW}
\int_{\Omega} \sigma_{ji} \delta \gamma_{ji} \textrm{d}\Omega + \int_{\Omega}  \rho \ddot{u}_i \delta u_{i} \textrm{d}\Omega  = \int_{\Gamma_{t}}t_{0i} \delta u_i \textrm{d}{\Gamma_{t}} + \int_\Omega b_{0i} \delta u_{i} \textrm{d}\Omega
\end{equation}  
where $\delta \gamma_{ij}= \delta u_{j,i} -e_{ijk} \delta \varphi_k $, $\delta u_i$ and $\delta \varphi_i$ are arbitrary strains, displacements, and rotations, respectively, satisfying compatibility equations with homogeneous essential boundary conditions. It must be observed here that the PVW in Equation \ref{PVW} is the weak formulation of both linear and angular momentum balances. Hence, the symmetry of the stress tensor and, consequently, the symmetry of the strain tensor are imposed in average sense. This is a significant difference compared to classical formulations for Cauchy continua in which the symmetry of the stress tensor is assumed ``a priori''. 

\section{Discontinuous Cell Method Approximation}

\subsection{Domain Discretization}
Let us consider a three-dimensional \emph{primal cell complex}, which, according to the customary terminology in algebraic topology \cite{tonti-2}, is a subdivision of the three-dimensional space $\mathds{R}^3$ through sets of vertices (0-cells), edges (1-cells), faces (2-cells), and volumes (3-cells). Next let us construct a \emph{dual cell complex} anchored to the primal. This can be achieved, for example, by associating a primal 3-cell with a dual 0-cell , a primal 2-cell with a dual 1-cell, etc. The primal/dual complex obtained through the Delaunay triangulation of a set of points and its associated Voronoi tessellation is a very popular choice in many fields of study for its ability to discretize complex geometry and it is adopted in this study. 

\begin{figure}[t]
\centering 
{\includegraphics[width=0.8\textwidth]{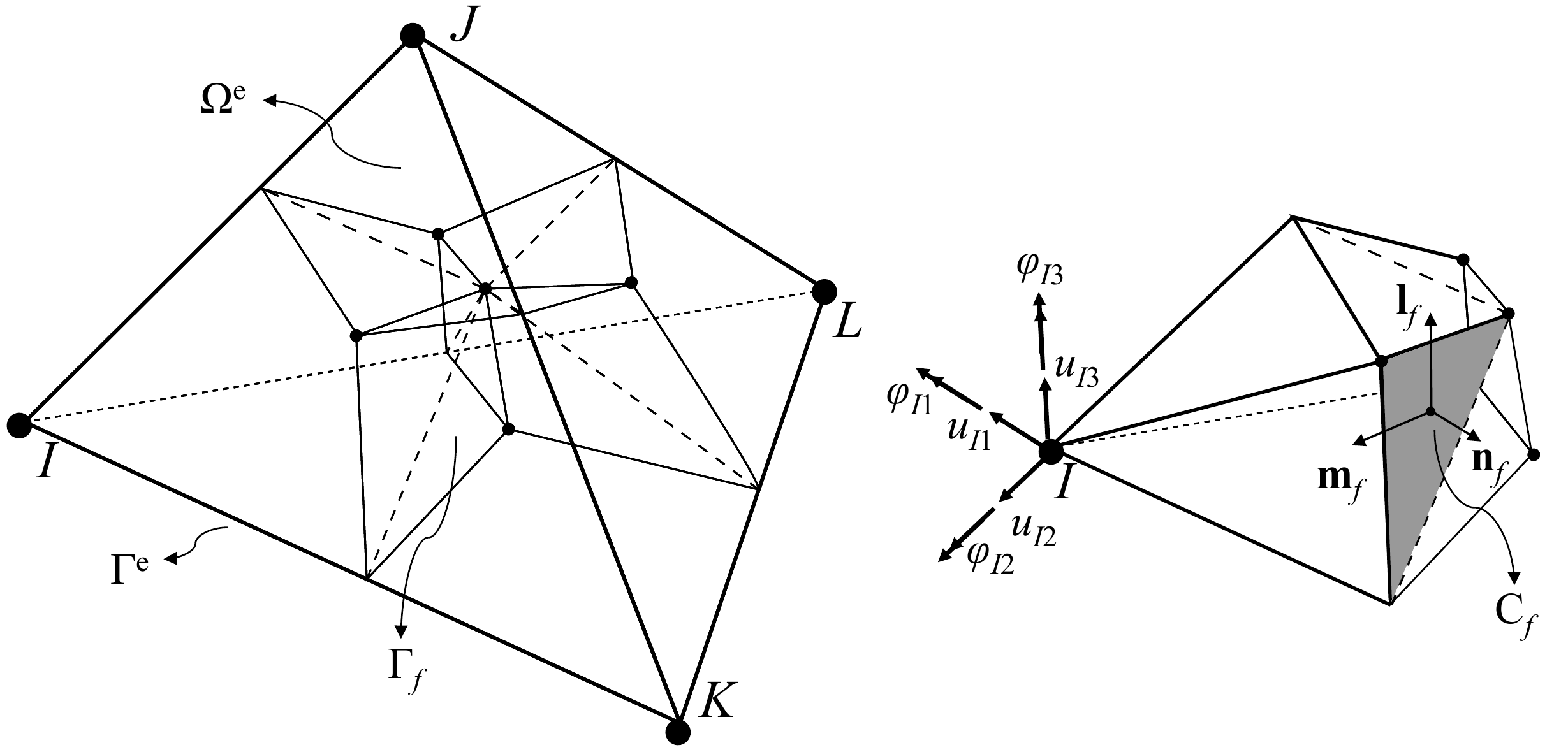}}
\caption{(a) Three dimensional Delaunay tetrahedralization and Voronoi tesselation. (b) Tetrahedron portion associated with node $I$. c) Voronoi cell. }
\label{DCMgeom3D}
\end{figure}

Let us consider a material domain $\Omega$ and discretize it into tetrahedral elements by using the centroidal Delaunay tetrahedralization, and the associated Voronoi tessellation which leads to a system of polyhedral cells \cite{Polymesh-paulino}. Figure \ref{DCMgeom3D}a, shows a typical tetrahedral element with the volume $\Omega^e$, external boundary $\Gamma^e$, and oriented surfaces $\Gamma^f$ located within the volume. The interior oriented surfaces $\Gamma^f$ are derived from the Voronoi tessellation and are hereinafter called ``facets". In 3D, the facets are triangular areas of contact between adjacent polyhedral cells. In the Voronoi tessellation procedure, the triangular facets $\Gamma^f$ are perpendicular to the element edges $\Gamma^e$, which is a crucial feature of DCM formulation as it will be shown later, and their geometry is such that one node of each facet is placed in the middle of the tetrahedral element edge, one is located on one of the triangular faces of the tetrahedral element, and one is located inside the tetrahedral element. As a result, each tetrahedral element contains twelve facets in a 3D setting, Figure \ref{DCMgeom3D}a. Figure \ref{DCMgeom3D}b illustrates a portion of the tetrahedral element associated with one of its four nodes $\alpha$ and the corresponding facets.  Combining such portions from all the tetrahedral elements connected to the same node, one obtains the corresponding Voronoi cell. Each node in the 3D DCM formulation has six degrees of freedom, three translational and three rotational, which are shown in Figure \ref{DCMgeom3D}b. The same figure depicts, for a generic facet, three unit vectors, one normal $\bold{n}_f$ and two tangential ones  $\bold{m}_f$ and $\bold{l}_f$, defining a local system of reference. In the rest of the paper, the facet index $f$ is dropped when possible to simplify notation.

\subsection{Discretized Kinematics}
The DCM approximation is based on the assumptions that displacement and rotation fields can be approximated by the rigid body kinematics of each Voronoi cell, that is 
\begin{equation}
u_i({\bf x}) = u_{Ii} + e_{ijk} \varphi_{Ij} (x_k - x_{Ik}) ;~~ \varphi _i ({\bf x}) = \varphi_{Ii} ~~~ \textrm{for} ~ {\bf x} \in \Omega^I
\label{eq:disp-approx}
\end{equation}
where $u_{Ii}$, $\varphi_{Ik}$ are displacements and rotations of  node $I$; $\Omega^I$ is the volume of the cell associated with node $I$.  Obviously with this approximating displacement and rotation functions, strain versus displacement/rotation relationships in Equations \ref{Eq:Compatibility-1}  cannot be enforced locally -- as typically done in displacement based finite element formulation.

Let us consider a generic node $I$ of spatial coordinates ${\bf x}_I$ and the adjacent nodes $J$ of spatial coordinates ${ \bf x}_J = { \bf x}_I+ \Bell$, where $\Bell$ is the vector connecting the two nodes. One can write $\Bell= \ell { \bf n}^\ell$ in which ${ \bf n}^\ell$ is a unit vector in the direction of $\Bell$. Note that, to simplify notation the two indices $I$ and $J$ are dropped when supposed to appear together. 
Without loss of generality, let us also assume that node $I$ is located on the negative side of the facets whereas node $J$ is located on the positive side of the associated facet oriented through its normal unit vector ${ \bf n}$. Moreover, it is useful to introduce the vector ${\bf d} = { \bf x}_{P} - { \bf x}_I$ connecting node $I$ to a generic point P on the facet. 
The displacement jump at point P reads
\begin{equation}
\hat{{\bf u}}= {\bf u}^+ - {\bf u}^- = {\bf u}_{J} - {\bf u}_I + \mb{\varphi}_{J} \times ({\bf d} - \Bell ) - \mb{\varphi}_I \times {\bf d} = \Delta {\bf u} - \mb{\varphi}_I \times \ell {\bf n}^\ell - \Delta \mb{\varphi} \times \ell ({\bf n}^\ell - \xi {\bf n}^d)\label{Eq:Strain1}
\end{equation}
where ${\bf u}^+$ and ${\bf u}^-$ are the values of displacements on the positive and negative side of the facet, respectively; $\Delta {\bf u} = {\bf u}_{J} - {\bf u}_I$; $\Delta \mb{\varphi} = \mb{\varphi}_{J} - \mb{\varphi}_I$; $\xi = d/\ell$; and $d$, ${\bf n}^d$ are magnitude and direction of vector ${\bf d}$. Equation \ref{Eq:Strain1} can be rewritten in tensorial notation as
\begin{equation}
\hat{u}_{i} =  \Delta u_{i} - \ell e_{ijk} \varphi_{j} n^\ell_{k} - \ell e_{ijk} \Delta \varphi_{j}  (n_{k}^\ell - \xi n_{k}^d) \label{Eq:Strain2}
\end{equation}

By expanding the displacement and the rotation fields in Taylor series around ${\bf x}_I$ and by truncating the displacement to the second order and the rotation to the first, one obtains
\begin{equation}
\Delta u_{i} = \ell u_{i,j} n^\ell_{j} + \frac{1}{2}\ell^2 u_{i,jk} n^\ell_{j} n^\ell_{k}; ~~~~ \Delta \varphi_{i} = \ell \varphi_{i,j} n^\ell_{j}
\end{equation}

By projecting the displacement jump $\hat{u}_{i}$ in the direction orthogonal to the facet and dividing it by the element edge length $\ell$, one can write

\begin{equation}
\begin{split}
\ell^{-1} \hat{u}_i n_i = \ell^{-1} n_i (\ell u_{i,j} n^\ell_j + \frac{1}{2}\ell^2 u_{i,jk} n^\ell_j n^\ell_k) - e_{ijk} n_i \varphi_j n^\ell_k - e_{ijk} n_i (\ell \varphi_{j,p} n^\ell_p)  (n_k^\ell - \xi n_k^d)= \\
= u_{j,i} n^\ell_i n_j + \frac{1}{2}\ell u_{j,ik} n_i^\ell n_j n^\ell_k - e_{ijk} \varphi_k n_i^\ell  n_j - \ell \varphi_{j,i} e_{pjk} n_p n^\ell_i  (n_k^\ell - \xi n_k^d) =\\
= (u_{j,i} - e_{ijk} \varphi_k) n^\ell_i n_j + \frac{1}{2}\ell u_{j,ik} n^\ell_i n_j n^\ell_k - \ell \varphi_{j,i} e_{pjk} n_p n^\ell_i  (n_k^\ell - \xi n_k^d) = \gamma_{ij}n^\ell_i n_j + \mathcal{O}(\ell)
\label{Eq:DisJumpN}
\end{split}
\end{equation}
where $\gamma_{ij}$ is the strain tensor. At convergence the discretization size tends to zero ($\ell \rightarrow 0$) and one can write $\ell^{-1} \hat{u}_i n_i =\gamma_{ij} n^\ell_i n_j$. Furthermore, if the dual complex adopted for the volume discretization is such that the facets are orthogonal to the element edges -- this condition is verified, for example, by the Delaunay-Voronoi complex -- then, ${ \bf n}^\ell \equiv { \bf n}$, and one has  $\ell^{-1} \hat{u}_i n_i =\gamma_{ij} n_i n_j$: the normal component of the displacement jump normalized with the element edge length represents the projection of the strain tensor onto the facet. Similarly, it can be shown that the components of the displacement jump tangential to the facets can be expressed as $\ell^{-1} \hat{u}_i m_i =\gamma_{ij} n_i m_j$ and $\ell^{-1} \hat{u}_i l_i =\gamma_{ij} n_i l_j$. 

Before moving forward, a few observations are in order. Since the facet are flat the unit vector $n_i$ is the same for any point belonging to a given facet and the projection of the strain tensor is uniform over each facet for a uniform strain field. The variation of the displacement jump over the facet is due to the curvature and it is an high order effect that can be neglected (see the last term in Equation \ref{Eq:DisJumpN}). Based on the previous observation one can conclude that the analysis of the interaction of two adjacent nodes can be based on the average displacement jump which, given the linear distribution of the jump, can be calculated as the displacement jump, $\bold{w}$, at the centroid of the facet C. This leads naturally to the following definition of ``facet strains'':
\begin{equation} \label{eq:eps-N}
\epsilon_{N} = \frac{ n_i } {\ell \Gamma}  \int_{\Gamma} \hat{u}_i \, \textrm{d}S = \frac{n_i w_{i}}{\ell} = n_i n_j \gamma_{ij}
\end{equation}
and 
\begin{equation}\label{eq:eps-ML}
\epsilon_{M} = \frac{m_i } {\ell \Gamma}  \int_{\Gamma} \hat{u}_i \, \textrm{d}S = \frac{m_i w_{i}}{\ell} = n_i m_j \gamma_{ij} ~~~~ \epsilon_{L} = \frac{l_i } {\ell \Gamma}  \int_{\Gamma} \hat{u}_i \, \textrm{d}S = \frac{l_i w_{i}}{\ell} = n_i l_j \gamma_{ij}
\end{equation}
Equations \ref{eq:eps-N} and \ref{eq:eps-ML} show that the ``facet strains'' correspond to the projection of the strain tensor onto the facet local system of reference.

Let us now consider a uniform hydrostatic stress/strain state in one element $\gamma_{ij} = \epsilon_V \delta_{ij}$ and $\sigma_{ij} = \sigma_V \delta_{ij}$. In this case the tractions on each facet must correspond to the volumetric stress and energetic consistency requires that $3 \Omega_e \sigma_{V} \epsilon_V = \sum_{ \mathcal{F}_e} \Gamma \ell \sigma_V \epsilon_N$ which gives 
\begin{equation} \label{eq:eps-V}
\epsilon_V=\frac{1}{3 \Omega_e} \sum_{ \mathcal{F}_e} \Gamma \ell \epsilon_N = \frac{1}{3 \Omega_e} \sum_{ \mathcal{F}_e} \Gamma n_i w_{i}
\end{equation}
where $\mathcal{F}_e$ is the set of facets belonging to one element. By using Equation \ref{eq:eps-V} and Equations \ref{eq:eps-N}, \ref{eq:eps-ML}, facet deviatoric strains can also be calculated as $\epsilon_D=\epsilon_N-\epsilon_V$.

By introducing Equation \ref{eq:disp-approx} into Equations \ref{eq:eps-N} and \ref{eq:eps-ML} one can also write 

\begin{equation} \label{eq:eps-N-new}
\ell \epsilon_{N} = n_i w_{i}= -n_{Ji} (u_{Ji} + e_{ijk} \varphi_{Jj} c_{Jk}) -n_{Ii} (u_{Ii} + e_{ijk} \varphi_{Ij} c_{Ik}) 
\end{equation}
\begin{equation} \label{eq:eps-M-new}
\ell \epsilon_{M} = m_i w_{i}= -m_{Ji} (u_{Ji} + e_{ijk} \varphi_{Jj} c_{Jk}) -m_{Ii} (u_{Ii} + e_{ijk} \varphi_{Ij} c_{Ik}) 
\end{equation}
\begin{equation} \label{eq:eps-L-new}
\ell \epsilon_{L} = l_i w_{i}= -l_{Ji} (u_{Ji} + e_{ijk} \varphi_{Jj} c_{Jk}) -l_{Ii} (u_{Ii} + e_{ijk} \varphi_{Ij} c_{Ik}) 
\end{equation}
where $\bold{n}_J=-\bold{n}$, $\bold{n}_I=\bold{n}$, $\bold{m}_J=-\bold{m}$, $\bold{m}_I=\bold{m}$, $\bold{l}_J=-\bold{l}$, $\bold{l}_I=\bold{l}$, and $\bold{c}_I$, $\bold{c}_J$ are vectors connecting the facet centroid C with nodes $I$, $J$, respectively.

\subsection{Discretized Equilibrium Equations}
For the adopted discretized kinematics in which all deformability is concentrated at the facets, the PVW in Equation \ref{PVW} can be rewritten as 

\begin{eqnarray} \label{eq:equi-weak-discrete-0}
\begin{split}
\sum_{\mathcal{F}} \Gamma \ell \left(  t_N \delta \epsilon_N + t_M \delta \epsilon_M + t_L \delta \epsilon_L  \right) + \int_{\Omega}  \rho \ddot{u}_i \delta u_{i} \textrm{d}\Omega= 
\int_{\Gamma_t}t_{0i} \delta u_i \textrm{d}{\Gamma} + \int_{\Omega} b_{0i} \delta u_{i} \textrm{d}\Omega 
\end{split}
\end{eqnarray}
where $\mathcal{F}$ is the set of all facets in the domain. 

By introducing Equations \ref{eq:disp-approx}, \ref{eq:eps-N-new}, \ref{eq:eps-M-new}, and \ref{eq:eps-L-new} into Equation \ref{eq:equi-weak-discrete-0} and considering displacement and traction continuity at the inter-element interfaces, one can write
\begin{eqnarray} \label{eq:equi-weak-discrete}
\begin{split}
\sum_I \left (  \sum_{\mathcal{F}_I} \Gamma t_{Ii} ( \delta u_{Ii} + e_{ijk}  \delta \varphi_{Ij} c_{Ik}) +   \int_{\Omega_I}   [b_{0i}-\rho \ddot{u}_{Ii} - \rho e_{imp} \ddot{\varphi}_{Im} (x_p-x_{Ip})] [ \delta u_{Ii} + e_{ijk} \delta \varphi_{Ij} (x_k-x_{Ik}) ] d\Omega \right)=0
\end{split}
\end{eqnarray}
where $t_{Ii} = t_N n_{Ii} +t_M m_{Ii}+t_L l_{Ii}$, $I$ is the generic Voronoi cell, $\Omega_I$ and $\mathcal{F}_I$ are the volume and the set of facets of the cell $I$, respectively. Note that the first term on the LHS of Equation \ref{eq:equi-weak-discrete} also includes the contribution of external tractions for cells located on the domain boundary.

Since Equation \ref{eq:equi-weak-discrete} must be satisfied for any virtual variation $ \delta u_{Ii}$ and $\delta \varphi_{Ik}$, it is equivalent to the following system of algebraic equations ($I=1, ..., N_c$, $N_c$= total number of Voronoi cells):

\begin{equation} \label{eq:cell-equilibrium-1}
\sum_{\mathcal{F}_I} \Gamma t_{Ii} + F_{Ii}-\mathcal{M}_{I} \ddot{u}_{Ii}- \mathcal{S}_{Iik} \ddot{\varphi}_{Ik}=0~~~\textrm{for}~~~i=1,2,3
\end{equation}
\begin{equation} \label{eq:cell-equilibrium-2}
\sum_{\mathcal{F}_I} \Gamma t_{Ii}e_{ijk} c_{Ij} + W_{Ik}-\mathcal{S}_{Iik} \ddot{u}_{Ii}- \mathcal{I}_{Ikp} \ddot{\varphi}_{Ip}=0~~~\textrm{for}~~~k=1,2,3
\end{equation}
where $F_{Ii}=\int_{\Omega_I} b_{0i} \textrm{d}\Omega$ = external force resultant, $\mathcal{M}_{I}=\int_{\Omega_I} \rho \textrm{d}\Omega$ = mass, $\mathcal{S}_{Iik} = \int_{\Omega_I} \rho e_{ijk} (x_j-x_{Ij}) \textrm{d}\Omega$ = first-order mass moments, $W_{Ik}=\int_{\Omega_I} b_{0i} e_{ijk} (x_j-x_{Ij}) \textrm{d}\Omega$ = external moment resultant, $\mathcal{I}_{Ikp}= \int_{\Omega_I} \rho e_{imp} (x_m-x_{Im}) e_{ijk} (x_j-x_{Ij}) \textrm{d}\Omega$ = second-order mass moments, of cell $I$. Equations \ref{eq:cell-equilibrium-1} and \ref{eq:cell-equilibrium-2} coincides with the force and moment equilibrium equations for each Voronoi cell. 

Note that $\mathcal{S}_{Iik}=0$ and $W_{Ik}=0$ (for uniform body force), if the vertex of the nodes of the Delaunay discretization coincide with the mass centroid of the Voronoi cells. This is the case for all the cells in the interior of the mesh if a centroidal Voronoi tessellation is adopted. Also, in general, $\mathcal{I}_{Ikp} \neq 0$ for $k\neq p$, and this leads to a non-diagonal mass matrix. A diagonalized mass matrix can be obtained simply by discarding the non-diagonal terms. 

\subsection{Discretized Constitutive Equations}
In the DCM framework, the constitutive equations are imposed at the facet level where the facet tractions need to be expressed as function of the facet strains. For elasticity, by projecting the tensorial constitutive equations reported in Equation \ref{eq:const} in the local system of reference of each facet, one has 

\begin{equation}\label{tN-const}
t_N=\sigma_{ij}n_i n_j= E_V \epsilon_V \delta_{ij}n_i n_j +E_D (\gamma_{ij}n_i n_j -  \epsilon_V \delta_{ij} n_i n_j) = E_V \epsilon_V + E_D \epsilon_D = E_D \frac{n_i w_{i}}{\ell} +  \frac{E_V-E_D}{3 \Omega_e} \sum_{ \mathcal{F}_e} \Gamma n_i w_{i}
\end{equation}
\begin{equation}\label{tM-const}
t_M=\sigma_{ij}n_i m_j= E_V \epsilon_V \delta_{ij}n_i m_j +E_D (\gamma_{ij}n_i m_j -  \epsilon_V \delta_{ij} n_i m_j) = E_D \epsilon_M = E_D \frac{m_i w_{i}}{\ell}
\end{equation}
\begin{equation}\label{tL-const}
t_L=\sigma_{ij}n_i l_j= E_V \epsilon_V \delta_{ij}n_i l_j +E_D (\gamma_{ij}n_i l_j -  \epsilon_V \delta_{ij} n_i l_j) = E_D \epsilon_L = E_D \frac{l_i w_{i}}{\ell}
\end{equation}

\section{Two-Dimensional Implementation}

\subsection{Three-Node Triangular Element}

\begin{figure}[t]
\centering 
{\includegraphics[width=\textwidth]{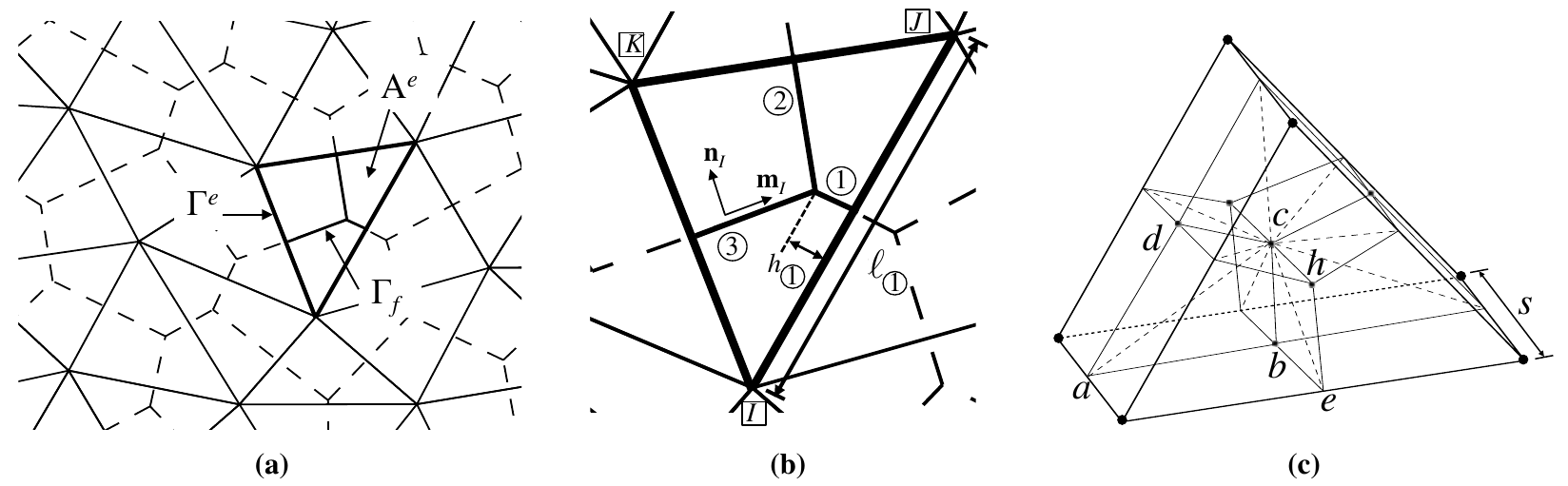}}
\caption{(a) Delaunay triangulation and Voronoi tesselation. (b) Facet geometry and associated edge length. (c) Triangular element geometry through the thickness.}
\label{DCMgeom}
\end{figure}

In order to pursue a two-dimensional implementation of DCM, let us consider a 2D Delaunay-Voronoi discretization as shown in Figures \ref{DCMgeom}a and b. A generic triangle can be considered as the triangular base of a prismatic volume as shown in Figure \ref{DCMgeom}c and characterized by 6 vertexes, 9 edges, 2 triangular faces, and 3 rectangular faces. By considering the Voronoi vertexes on the two parallel triangular faces and face/edge points located at mid-thickness (see, e.g., points $a$ and $d$) a complete tessellation of the volume in six sub-volumes, one per vertex, can be obtained by triangular facets.  Of these facets, $N_f=12$ are orthogonal to the triangular faces, set $\mathcal{F}_f$ (see, e.g., the one connecting points $c$, $e$, $h$ in Figure \ref{DCMgeom}c) and $N_o=6$ are parallel to the triangular faces, set $\mathcal{F}_o$,  (see, e.g., the one connecting points $c$, $d$, $a$ in Figure \ref{DCMgeom}c). 

One can write 
\begin{equation} \label{eq:eps-V-2D}
\epsilon_V=\frac{1}{3 \Omega_e} \left( \sum_{ \mathcal{F}_f} \Gamma \ell \epsilon_N + \sum_{ \mathcal{F}_o} \Gamma s \epsilon_N  \right)
\end{equation}
where $s$ is the out-of-plane thickness.

For plane strain conditions $\epsilon_N=0$ for the facet set $\mathcal{F}_o$ and simply the second term in Equation \ref{eq:eps-V-2D} is zero. For plane stress, instead, $t_N=E_V \epsilon_V+E_D(\epsilon_N-\epsilon_V)=0$ for the facet set $\mathcal{F}_o$. Therefore, for the facet set $\mathcal{F}_o$, $\epsilon_N=(E_D-E_V) \epsilon_V / E_D= -3\nu / (1-2 \nu) \epsilon_V $ and , also, $\sum_{ \mathcal{F}_o} \Gamma s = \Omega_e$. Using these relations, one has $\sum_{ \mathcal{F}_o} \Gamma s \epsilon_N = -3\nu/(1-2\nu)\Omega_e\epsilon_V$. Substituting this relation in Equation \ref{eq:eps-V-2D} leads to $\epsilon_V = (1-2\nu)/(3\Omega_e(1-\nu)) \left( \sum_{ \mathcal{F}_f} \Gamma \ell \epsilon_N   \right)$. In addition the facet set $\mathcal{F}_f$ is composed by 3 sets of 4 planar triangular facets. For each set, strains and tractions are the same on the 4 facets because the response is uniform through the thickness. Consequently the 4 facets can be grouped into one rectangular facet of area $sh$ where $h$ is the in length of the facet (see Figure \ref{DCMgeom}b)

By taking everything into account the volumetric strain for 2D problem can be written as 
\begin{equation} \label{eq:eps-V-2D-1}
\epsilon_V=\frac{1}{3 A^e \alpha} \sum_{f=1}^3 \ell_f h_f \epsilon_{fN} = \frac{1}{3 A^e \alpha} \sum_{f=1}^3 h_f n_{fi} w_{fi}
\end{equation}
where $A^e$ is the area of the triangular element and $\alpha=(1-\nu)/(1-2\nu)$ for plane stress and $\alpha=1$ for plane strain.  

The triangular DCM element has 9 degrees of freedom, two displacements and one rotation for each node, which can be collected in one vector $\bold{Q}^T=[u_{I1}~u_{I2}~\varphi_{I3}~u_{J1}~u_{J2}~\varphi_{J3}~u_{K1}~u_{K2}~\varphi_{K3}]$ in which $I$, $J$, and $K$ are the element node indexes, and 1, 2, and 3 represent the three Cartesian coordinate axes. By using Equations \ref{eq:eps-N-new} to \ref{eq:eps-L-new}, one can write $\epsilon_{fN}=\bold{N}_f \bold{Q}$, $\epsilon_{fM}=\bold{M}_f \bold{Q}$, $\epsilon_{V}=\bold{V} \bold{Q}$, $\epsilon_{fD}=\bold{D}_f \bold{Q}$ where
\begin{equation} \label{eq:matrix-N1}
\bold{N}_{\circled{1}} = \ell^{-1}_{\circled{1}} \begin{bmatrix}
-n_{I1} & -n_{I2} & n_{I1}c_{I2}-n_{I2}c_{I1} & n_{I1} & n_{I2} & -n_{I1}c_{J2}+n_{I2}c_{J1} & 0 & 0 & 0
\end{bmatrix}
\end{equation}
\begin{equation} \label{eq:matrix-M1}
\bold{M}_{\circled{1}} = \ell^{-1}_{\circled{1}} \begin{bmatrix}
-m_{I1} & -m_{I2} & m_{I1}c_{I2}-m_{I2}c_{I1} & m_{I1} & m_{I2} & -m_{I1}c_{J2}+m_{I2}c_{J1} & 0 & 0 & 0
\end{bmatrix}
\end{equation}
\begin{equation} \label{eq:matrix-N2}
\bold{N}_{\circled{2}} = \ell^{-1}_{\circled{2}} \begin{bmatrix}
0 & 0 & 0 & -n_{J1} & -n_{J2} & n_{J1}c_{J2}-n_{J2}c_{J1} & n_{J1} & n_{J2} & -n_{J1}c_{K2}+n_{J2}c_{K1}
\end{bmatrix}
\end{equation}
\begin{equation} \label{eq:matrix-M2}
\bold{M}_{\circled{2}} = \ell^{-1}_{\circled{2}} \begin{bmatrix}
0 & 0 & 0 & -m_{J1} & -m_{J2} & m_{J1}c_{J2}-m_{J2}c_{J1} & m_{J1} & m_{J2} & -m_{J1}c_{K2}+m_{J2}c_{K1}
\end{bmatrix}
\end{equation}
\begin{equation} \label{eq:matrix-N3}
\bold{N}_{\circled{3}} = \ell^{-1}_{\circled{3}} \begin{bmatrix}
-n_{I1} & -n_{I2} & n_{I1}c_{I2}-n_{I2}c_{I1} & 0 & 0 & 0 & n_{I1} & n_{I2} & -n_{I1}c_{K2}+n_{I2}c_{K1}
\end{bmatrix}
\end{equation}
\begin{equation} \label{eq:matrix-M3}
\bold{M}_{\circled{3}} = \ell^{-1}_{\circled{3}} \begin{bmatrix}
-m_{I1} & -m_{I2} & m_{I1}c_{I2}-m_{I2}c_{I1} & 0 & 0 & 0 & m_{I1} & m_{I2} & -m_{I1}c_{K2}+m_{I2}c_{K1}
\end{bmatrix}
\end{equation}
and $\bold{V} = (3 A_f \alpha)^{-1}\sum_{f=1}^3 \ell_f h_f \bold{N}_f$, $\bold{D}_f=\bold{N}_f-\bold{V}$

The elastic stiffness matrix for a single DCM triangle can be computed similarly to classical FEM elements by writing the internal energy as function of the nodal degrees of freedom. By using Equations \ref{tN-const} and \ref{tM-const}, one has 
\begin{equation} \label{IntEneVar-0}
\mathcal{U} = \frac{1}{2} s \ell_f h_f \left( t_{Nf} \epsilon_{Nf}  + t_{Mf} \epsilon_{Mf}\right)= \frac{1}{2} s \ell_f h_f \left[ E_V\epsilon^2_{Vf} +E_D\epsilon_{Df}^2+E_D\epsilon_{Mf}^2 + (E_V+E_D)\epsilon_{Vf}\epsilon_{Df}\right] =\frac{1}{2}\bold{Q}^T\bold{K} \bold{Q}
\end{equation}
where 
\begin{equation} \label{eq:K}
\bold{K} =  2s A^e E_V \bold{V}^T \bold{V} + sh_f \ell_f \left[ E_D \bold{D}_f^T \bold{D}_f + 0.5(E_V+E_D)(\bold{V}^T \bold{D}_f + \bold{D}_f \bold{V}^T) + E_D \bold{M}_f^T \bold{M}_f \right]
\end{equation}
and summation rule applies for repeated index $f$. In the previous equation the relation $\sum_f h_f \ell_f = 2A^e$ was used.

Similarly the mass matrix of the element can be obtained by writing the kinetic energy as function of the nodal velocities. For uniform density, one has 
\begin{equation} \label{node-kinetic}
\mathcal{K} = \frac{1}{2}  \sum_I\int_{A_{I}}  \rho (\dot{u}_{1}^2 + \dot{u}_{2}^2 ) s \textrm{d}A =  \frac{1}{2} s \rho \sum_I\int_{A_{I}} [\dot{u}^2_{I1} + \tilde{y}_I^2\dot{\varphi}^2_{I3}-2\tilde{y}_I\dot{u}_{I1}\dot{\varphi}_{I3}+ \dot{u}^2_{I2} + \tilde{x}^2_I \dot{\varphi}_{I3}^2+2\tilde{x}_I \dot{u}_{I2}\dot{\varphi}_{I3}] \textrm{d}A = \frac{1}{2}\dot{\bold{Q}}^T\bold{M} \dot{\bold{Q}}
\end{equation}
where $\tilde{x}_I=x - x_{I}$, $\tilde{y}_I=y - y_{I}$
\begin{equation} \label{element-mass}
\bold{M} = 
\begin{bmatrix}
\bold{M}_1 & \bold{0} & \bold{0} \\
\bold{0} & \bold{M}_2 & \bold{0} \\
\bold{0} & \bold{0} & \bold{M}_3 \\
\end{bmatrix}
\hspace{0.25 in}
\bold{M}_I =\rho s 
\begin{bmatrix}
A_I & 0 & -S_{yI} \\
0 & A_I &  S_{xI} \\
-S_{yI} & S_{xI} & I_I \\
\end{bmatrix}
\end{equation}
and $A_I=\int_{A_{I}} \textrm{d}A$, $I_I=\int_{A_{I}} (\tilde{x}^2_I+\tilde{y}^2_I )  \textrm{d}A$, $S_{yI}=\int_{A_{I}} \tilde{y}_I \textrm{d}A$, $S_{xI}=\int_{A_{I}} \tilde{x}_I \textrm{d}A$

\begin{figure}[t]
\centering 
{\includegraphics[width=0.85\textwidth]{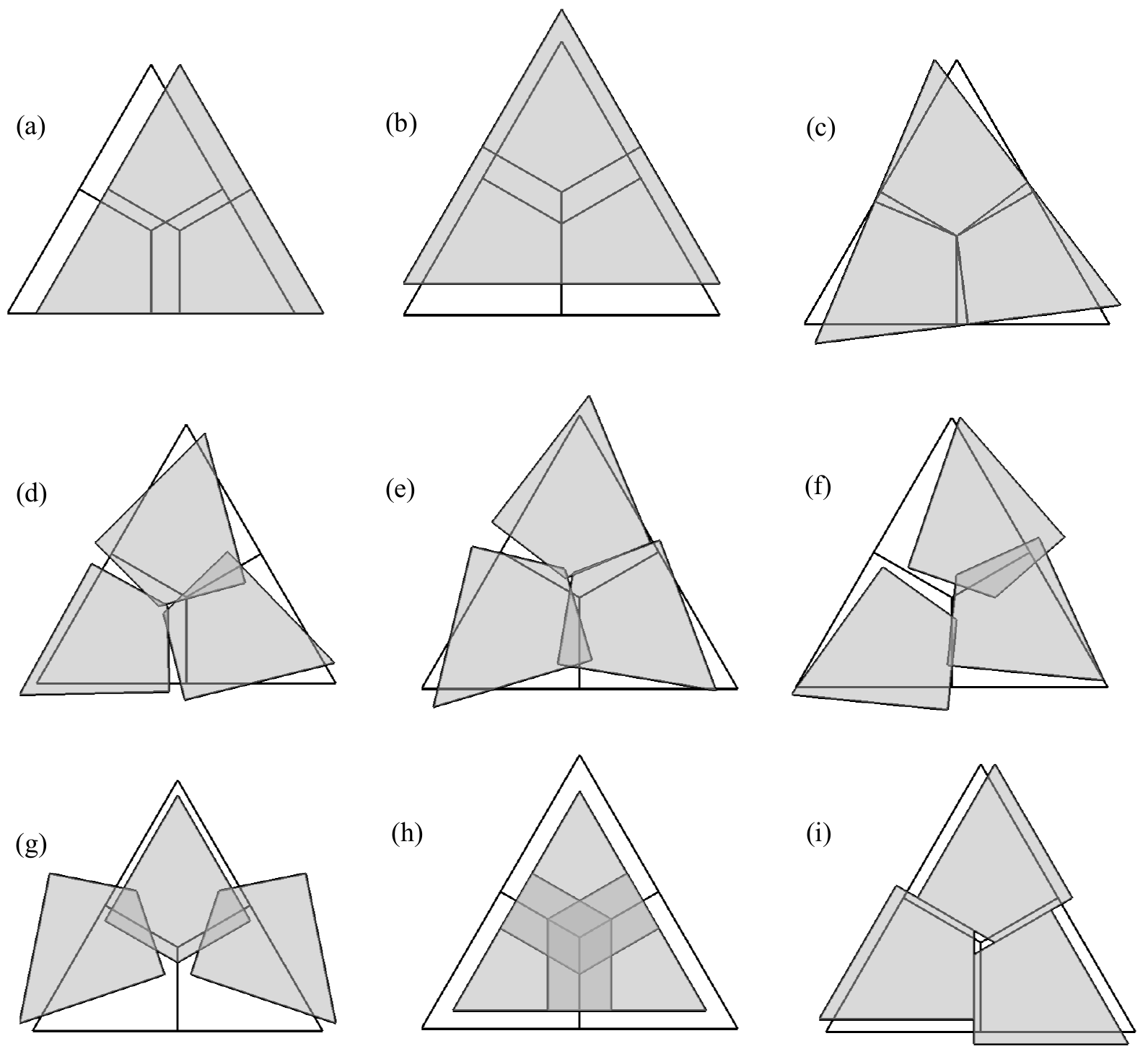}}
\caption{Deformation modes of an equilateral DCM triangular element: (a) to (c) rigid body modes, (d) to (f) bending modes, (g) uniaxial tension mode, (h) volumetric mode, and (i) pure shear mode.}
\label{def-modes}
\end{figure}

To study the deformation properties of the formulated triangular elements, one can consider the generalized eigenvalue problem $\bold K_e \boldsymbol \Phi_I  = \lambda_I \boldsymbol \Phi_I$ where $\lambda_I$ and $\boldsymbol \phi_I$ are the $I$th eigenvalue and the corresponding eigenvector of the element stiffness matrix, respectively. The triangular element has nine eigenvalues which correspond to the element nine degrees of freedom, three of which must be equal to zero to represent the three possible rigid body deformation modes. The other eigenvalues must be positive to ensure positive definiteness of the stiffness matrix. The deformation modes for an equilateral triangular element are plotted in Figure \ref{def-modes}: Figures \ref{def-modes} (a), (b), and (c) correspond to the rigid body deformation modes with zero eigenvalues; Figures \ref{def-modes} (d), (e), and (f) are bending modes of deformation; uniaxial deformation mode is depicted in Figures \ref{def-modes} (g); and volumetric and pure shear deformation modes are plotted in Figures \ref{def-modes} (h) and (i), respectively.

\begin{figure}[t]
\centering 
{\includegraphics[width=0.9\textwidth]{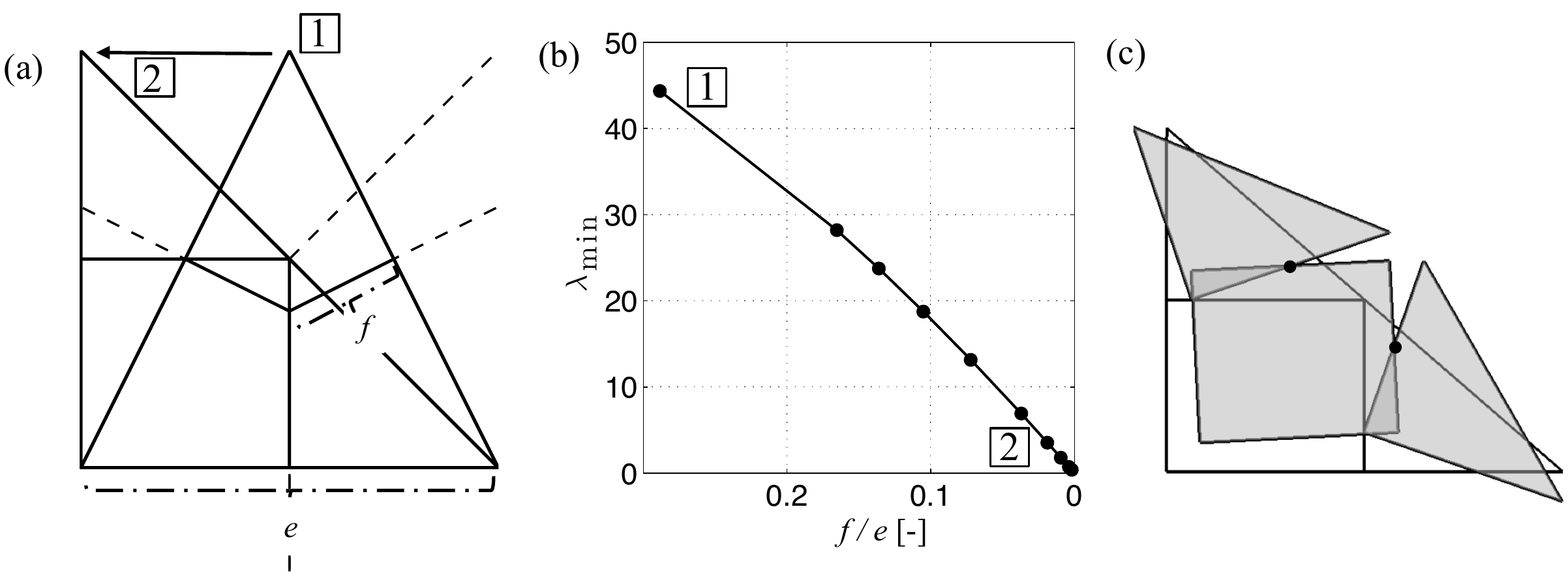}}
\caption{Variation of the minimum eigen value of the stiffness matrix with respect to the element geometry change.}
\label{spectralAnalysis}
\end{figure}

In Figure \ref{spectralAnalysis}a, a triangular element is considered, and the eigenvalue problem is solved for the different position of the top node from 1 to 2. For the final configuration 2, the length of the internal facet $f$ is zero. The minimum positive eigenvalue $\lambda_{min}$ is plotted versus facet to edge length ratio, $f/e$, of the triangular element in Figure \ref{spectralAnalysis}b. One can see that as the facet length tends to zero and the element becomes a right triangle, $\lambda_{min}$ tends to zero, which means that the element stiffness matrix becomes singular. The zero eigenvalue is associated with a zero energy mode of deformation, which is plotted in Figure \ref{spectralAnalysis}c. Since normal and tangential components of displacement jump vector at centroid of the facets are zero, normal and tangential strains on both element facets are equal to zero. Therefore, this element deformation mode is of zero energy. For $f/e < 0$, $\lambda_{min} < 0$ and the stiffness matrix is not positive definite. Therefore, during the domain discretization procedure, it must be considered that all element facets length must be positive. In other words, they must be completely placed inside the element, which, in turn, requires that all angles of the triangle to be smaller than $\pi/2$.

\subsection{Four-Node Quadrilateral Elements}
\begin{figure}[t]
\centering 
{\includegraphics[width=\textwidth]{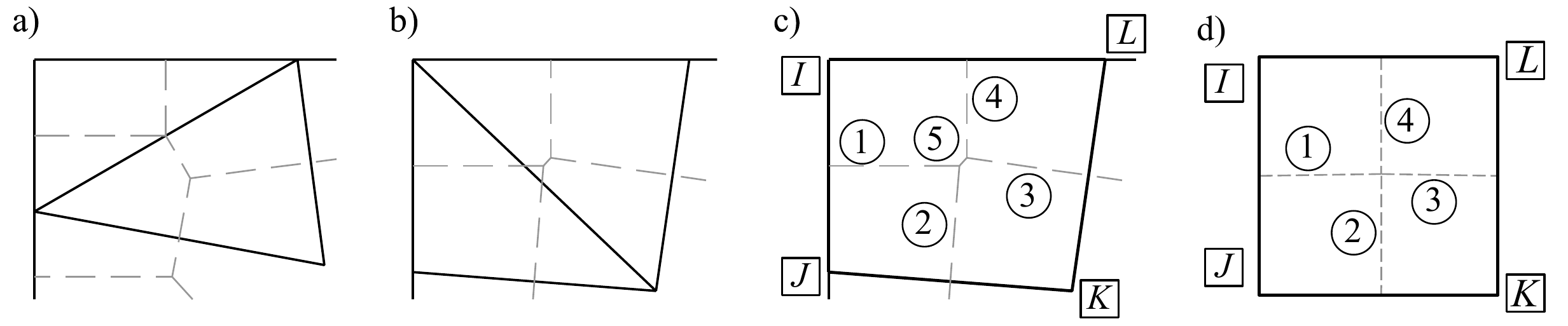}}
\caption{(a) Right triangle. (b) Obtuse triangle. (c) Typical quadrilateral element obtained from combining two adjacent triangular elements. (d) Quadrilateral element obtained from combining two adjacent right triangular elements.}
\label{34elems}
\end{figure}

It is not always possible to obtain triangular meshes that satisfy the shape requirements discussed in the previous section. Figure \ref{34elems}a shows, for example, a situation where the triangulation produced a right triangle at the right-angled exterior corner of a structural domain. In this situation, the Voronoi tessellation produces three orthogonal bisectors that intersect at a point located exactly on the hypotenuse of the right triangle. This results in a facet with zero length, which would lead to a zero eigenvalue of the stiffness matrix. Figure \ref{34elems}b shows a more severe situation where the generation of nodes and the resulting triangulation produces an obtuse triangle. In this situation, the three orthogonal bisectors intersect at a point located outside of the obtuse triangle, which would result in a negative facet area and a negative eigenvalue of the stiffness matrix. In order to overcome the computational problems implied by this situations, it is possible to combine a problematic triangle with the neighboring triangle into a four-node, five-facet quadrilateral element, Figure \ref{34elems}c. The interior, or fifth facet is the orthogonal bisector of a straight line between two nodes at opposite corners of the quadrilateral element. The translation and rotation of the two nodes labeled as $I$ and $K$ in Figure \ref{34elems}c, are the degrees of freedom that produce the displacement jump for the fifth facet. The straight line distance between Nodes $I$ and $K$ provides the edge length $\ell_{\circled{5}}$ required for the calculation of facet strains. The volumetric strain, constant inside the quadrilateral element, is still calculated by Equation \ref{eq:eps-V-2D-1} where the contribution of all five facets are taken into account. Also, a four node, rectangular element with four facets is generated if the two adjacent triangles are both right, see Figure \ref{34elems}d. For the generic quadrilateral element, Equations \ref{eq:matrix-N1} to \ref{eq:matrix-M3} must be substituted by the following equations

\begin{equation} \label{4eq:matrix-N1}
\bold{N}_{\circled{1}} = \ell^{-1}_{\circled{1}} \begin{bmatrix}
-n_{I1} & -n_{I2} & n_{I1}c_{I2}-n_{I2}c_{I1} & n_{I1} & n_{I2} & -n_{I1}c_{J2}+n_{I2}c_{J1} & 0 & 0 & 0 & 0 & 0 & 0
\end{bmatrix}
\end{equation}
\begin{equation} \label{4eq:matrix-M1}
\bold{M}_{\circled{1}} = \ell^{-1}_{\circled{1}} \begin{bmatrix}
-m_{I1} & -m_{I2} & m_{I1}c_{I2}-m_{I2}c_{I1} & m_{I1} & m_{I2} & -m_{I1}c_{J2}+m_{I2}c_{J1} & 0 & 0 & 0 & 0 & 0 & 0
\end{bmatrix}
\end{equation}
\begin{equation} \label{4eq:matrix-N2}
\bold{N}_{\circled{2}} = \ell^{-1}_{\circled{2}} \begin{bmatrix}
0 & 0 & 0 & -n_{J1} & -n_{J2} & n_{J1}c_{J2}-n_{J2}c_{J1} & n_{J1} & n_{J2} & -n_{J1}c_{K2}+n_{J2}c_{K1} & 0 & 0 & 0
\end{bmatrix}
\end{equation}
\begin{equation} \label{4eq:matrix-M2}
\bold{M}_{\circled{2}} = \ell^{-1}_{\circled{2}} \begin{bmatrix}
0 & 0 & 0 & -m_{J1} & -m_{J2} & m_{J1}c_{J2}-m_{J2}c_{J1} & m_{J1} & m_{J2} & -m_{J1}c_{K2}+m_{J2}c_{K1} & 0 & 0 & 0
\end{bmatrix}
\end{equation}
\begin{equation} \label{4eq:matrix-N3}
\bold{N}_{\circled{3}} = \ell^{-1}_{\circled{3}} \begin{bmatrix}
0 & 0 & 0 & 0 & 0 & 0 & -n_{K1} & -n_{K2} & n_{K1}c_{K2}-n_{K2}c_{K1} & n_{K1} & n_{K2} & -n_{K1}c_{L2}+n_{K2}c_{L1}
\end{bmatrix}
\end{equation}
\begin{equation} \label{4eq:matrix-M3}
\bold{M}_{\circled{3}} = \ell^{-1}_{\circled{3}} \begin{bmatrix}
0 & 0 & 0 & 0 & 0 & 0 & -m_{K1} & -m_{K2} & m_{K1}c_{K2}-m_{K2}c_{K1} & m_{K1} & m_{K2} & -m_{K1}c_{L2}+m_{K2}c_{L1}
\end{bmatrix}
\end{equation}
\begin{equation} \label{4eq:matrix-N4}
\bold{N}_{\circled{4}} = \ell^{-1}_{\circled{4}} \begin{bmatrix}
-n_{I1} & -n_{I2} & n_{I1}c_{I2}-n_{I2}c_{I1} & 0 & 0 & 0 & 0 & 0 & 0 & n_{I1} & n_{I2} & -n_{I1}c_{L2}+n_{I2}c_{L1}
\end{bmatrix}
\end{equation}
\begin{equation} \label{4eq:matrix-M4}
\bold{M}_{\circled{4}} = \ell^{-1}_{\circled{4}} \begin{bmatrix}
-m_{I1} & -m_{I2} & m_{I1}c_{I2}-m_{I2}c_{I1} & 0 & 0 & 0 & 0 & 0 & 0 & m_{I1} & m_{I2} & -m_{I1}c_{L2}+m_{I2}c_{L1}
\end{bmatrix}
\end{equation}
\begin{equation} \label{4eq:matrix-N5}
\bold{N}_{\circled{5}} = \ell^{-1}_{\circled{5}} \begin{bmatrix}
-n_{I1} & -n_{I2} & n_{I1}c_{I2}-n_{I2}c_{I1} & 0 & 0 & 0 & n_{I1} & n_{I2} & -n_{I1}c_{K2}+n_{I2}c_{K1} & 0 & 0 & 0
\end{bmatrix}
\end{equation}
\begin{equation} \label{4eq:matrix-M5}
\bold{M}_{\circled{5}} = \ell^{-1}_{\circled{5}} \begin{bmatrix}
-m_{I1} & -m_{I2} & m_{I1}c_{I2}-m_{I2}c_{I1} & 0 & 0 & 0 & m_{I1} & m_{I2} & -m_{I1}c_{K2}+m_{I2}c_{K1} & 0 & 0 & 0
\end{bmatrix}
\end{equation}

\section{Elastic Analysis Results}
\begin{figure}[t]
\centering 
{\includegraphics[width=\textwidth]{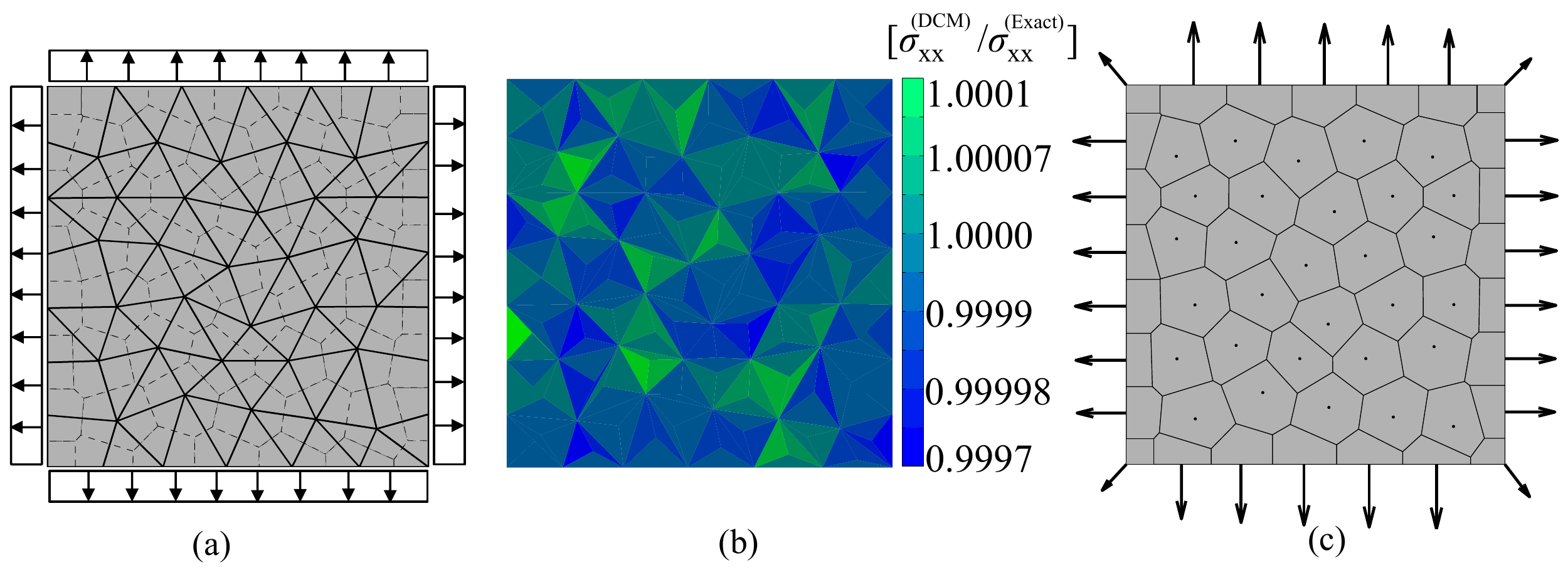}}
\caption{(a) Unit square used in the patch test subjected to uniform bi-axial strain field. (b) Contour of $\sigma_{xx}$ obtained from the patch test normalized with the exact value. (c) Vector representation of nodal forces obtained from the patch test.}
\label{PatchTest}
\end{figure}

\subsection{Patch Test and Facet Tensors}
Numerical experiments carried out in this section show that the 2D DCM triangle passes the patch test and is able to reproduce exactly uniform strain and stress fields. Acknowledging that the following observations apply equally to stress or strain, the most basic result regarding a uniform field is that the normal and tangential stresses calculated by DCM for a facet with certain orientation correspond to the tractions calculated by projecting the stress tensor onto the facet orientation. Conversely, facet's orientation, normal and tangential (shear) stresses along with the calculated value of the volumetric stress can be used to determine the facet overall stress tensor with respect to the global system of reference (see Appendix \ref{facet strain tensor}). To do the patch test, an elastic square specimen discretized by DCM is subjected to uniform $\varepsilon_{xx} = \varepsilon_{yy} = 0.1$ as shown in Figure \ref{PatchTest}a. $E = 1000$ MPa and $\nu$ = 0.25 are used as material properties. DCM analysis is performed, and the results are used to calculated the stress tensor for all of the facets. The exact uniform stress tensor due to the applied uniform strain tensor can be calculated by elasticity equations. In Figure \ref{PatchTest}b, $\sigma_{xx}^{\text{DCM}}/\sigma_{xx}^{\text{Exact}} =$ ratio of the $\sigma_{xx}$ calculated by DCM for each facet to the one obtained from elasticity equations, is plotted for all facets. For each facet and the corresponding element, the portion of the element area associated to that facet is colored according to the $\sigma_{xx}^{\text{DCM}}/\sigma_{xx}^{\text{Exact}}$ value. One can see that DCM successfully generated the uniform stress field which matches well with the elasticity results. This contour is the same for the $\sigma_{yy}$ component of the stress tensor. In addition, the resultant force vector on each node is plotted in Figure \ref{PatchTest}c. It is clear that the force vector is negligible for the nodes inside the specimen, while its distribution on the specimen surfaces correspond to the uniform stress and strain fields. 

\subsection{Convergence Study on Cantilever Beam} \label{cantilever}
\begin{figure}[h!]
\centering 
{\includegraphics[width=0.7\textwidth]{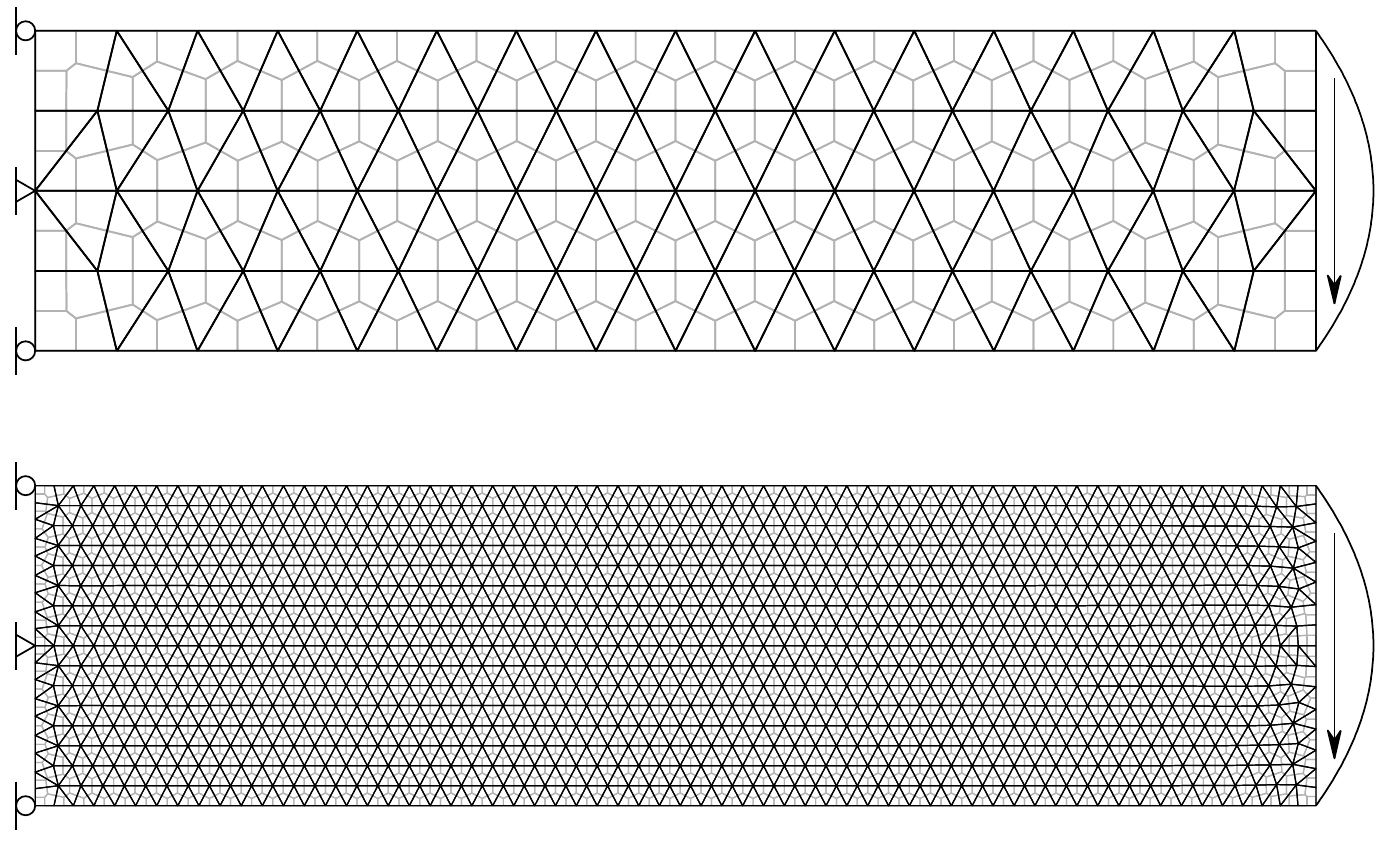}}
\caption{Coarsest (top) and finest (bottom) mesh used for the convergence study.}
\label{CanteliverBeam}
\end{figure}

In order to study the convergence of the present method to the exact solution for a non-uniform strain field, a classical cantilever beam test [19] was simulated. The rectangular domain, shown in Figure \ref{CanteliverBeam}, is characterized by a length-to-depth ratio of 4. The traction boundary conditions are the classic stress distributions of simple bending. Figure \ref{CanteliverBeam} shows a parabolically varying shear at the cantilever tip. At the fixed end, only the displacement boundary conditions are shown for clarity. However, an equal but opposite parabolic shear was applied at the fixed end, as well as a linearly varying normal stress based on the non-zero bending moment at that location. The exact solution for the displacement field is provided in Hughes [19] which assumes linear isotropic elasticity.

\begin{figure}[t]
\centering 
{\includegraphics[width=0.8\textwidth]{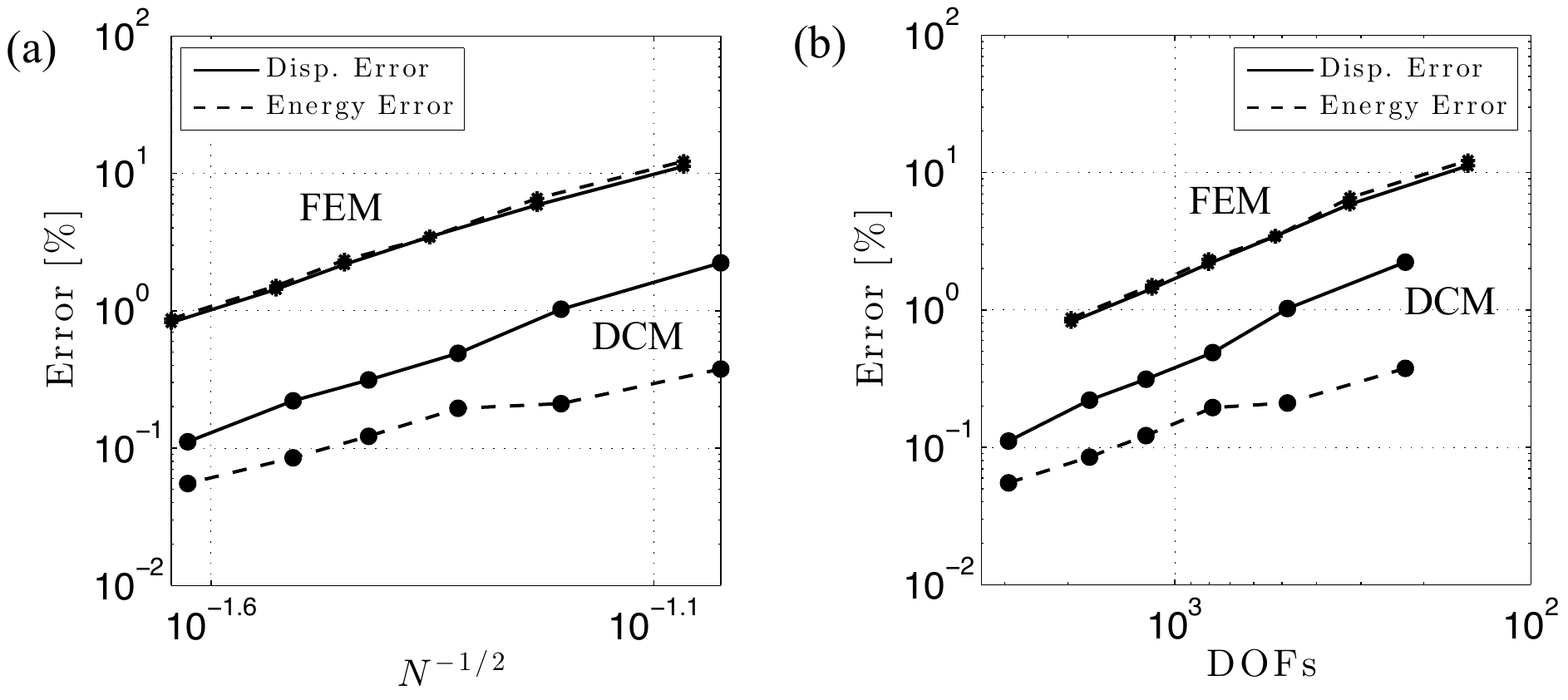}}
\caption{Convergence Study of the cantilever beam regarding beam tip displacement error and total strain energy error.}
\label{PeakStat}
\end{figure}

Six different meshes at various levels of refinement are considered. Figure \ref{CanteliverBeam} shows the coarsest, with 128 elements, and the finest, with 1790 elements. For comparison, the same numerical simulations were performed using the standard constant strain triangle (CST) finite element. All the computations were carried out under plane strain conditions, with a Poisson's ratio of 0.3. Figure \ref{PeakStat} presents the results of the convergence study. The relative error between the numerical calculation and the exact solution is plotted as a function of the inverse of the square root of the number of elements ($N^{−1/2}$), which is proportional to the characteristic element size. The results for the total elastic energy and the tip displacement are shown in Figures \ref{PeakStat}a for both DCM and Constant Strain Triangle (CST) finite element.

In the log-log plots of Figures \ref{PeakStat}a, the slope of the line segments provide a measure of the convergence rate. For the strain energy, the average convergence rates for the DCM and CST, respectively, are 1.62 and 1.99, and for the tip displacement they are 2.1 and 2.02. The theoretical convergence rate for the CST is 2 for both strain energy and tip deflection. Although the convergence rates are comparable, the DCM outperforms CST in terms of accuracy. The CST error in both strain energy and tip deflection is one order of magnitude higher than the DCM error. In terms of energy, for example, the DCM error ranges from 0.38$\%$ to 0.06$\%$ (coarsest to finest mesh), whereas the CST error ranges from 10$\%$ to 0.8$\%$. It must be mentioned, however, that each node in the DCM has one degree of freedom, the rotation, more than its counterpart in the CST. This additional degree of freedom results in higher computational cost for DCM compared to classical FEM. The tip displacement and the strain energy errors are plotted versus total number of degrees of freedom for both DCM and FEM simulations in Figure \ref{PeakStat}b in log-log axes. One can see that as the total number of DOFs increases, the error values decrease (DOFs axis is reversed). It can bee seen that for approximately equal number of elements, total number of DOFs for DCM is higher than CST. However, the accuracy remains higher than CST for the same number of DOFs.


\section{Cohesive Fracture Propagation}
The convergence study presented in the previous section demonstrates that the DCM performs very well in the elastic regime. However, the most attractive feature of this method is the ability of easily accommodating the displacement discontinuity associated with fracture without suffering from the typical shortcomings of the classical finite element method, the limitations of typical particle models, or the complexity and the high computational cost of advanced finite element formulations. In this section a simple isotropic damage model is introduced in the DCM framework in order to simulate the initiation and propagation of quasi-brittle fracture.

\subsection{Formulation}
According to the classical damage mechanics and the DCM formulation for elasticity presented above, in a damaged material the facet tractions $t_{fN}$ and $t_{fM}$ can be calculated as:
\begin{equation} \label{NTConstDamag}
t_{fN} = (1-D_f) \bigg[ \frac{E_V-E_D}{3\alpha_v\Omega^e}\sum_f A_{f} w_{fN} + \frac{E_D w_{fN}}{\ell_f} \bigg]; ~~~~~ t_{fM} = (1-D_f) \bigg[ \frac{E_D w_{fM}}{\ell_f} \bigg]
\end{equation}

where $D_f$ is the damage parameter related to the facet $f$. The evolution of the damage parameter is assumed to be governed by a history variable, the facet maximum effective strain, $\varepsilon_{f}^{max}$, characterizing the overall amount of straining that the material has been subject to during prior loading:
\begin{equation} \label{DamagPara}
D_f = 1 - \frac{\varepsilon_t}{\varepsilon_{f}^{max}} \text{exp} \bigg[\frac{-<\varepsilon_{f}^{max} - \varepsilon_t>}{\varepsilon_{fF}} \bigg]
\end{equation}
where $\left\langle{x}\right\rangle  = \text{max}(0,x)$, $\varepsilon_t$ is a material parameter representing the strain limit which governs the onset of damage, and $\varepsilon_{fF}$ governs the damage evolution rate. The maximum effective strain $\varepsilon_{f}^{max}$ that is used for each facet at each computational step is equal to the maximum principal strain $\varepsilon_{fI}$ that the facet has experienced through the loading process. This value is compared to $\varepsilon_t$ to distinguish facet elastic behavior from the the nonlinear case. $\varepsilon_{fI}$ is the maximum eigenvalue of the facet strain tensor whose components can be derived in terms of facet normal $\varepsilon_{fN}$, tangential $\varepsilon_{fM}$, and volumetric $\varepsilon_{V}$ strains as discussed in Section \ref{cantilever} and presented in details in Appendix \ref{facet strain tensor}.

In order to ensure convergence upon mesh refinement and to avoid spurious mesh sensitivity, one can define
\begin{equation} \label{effmaxstr}
\varepsilon_{fF} = \frac{\varepsilon_t}{2}\bigg( \frac{\ell_t}{\ell_f} - 1 \bigg)
\end{equation}
where $\ell_t = 2EG_t/\sigma_t^2$ is Hillerborg's characteristic length, which is assumed to be a material parameter. $\sigma_t$ and $G_t$ are the elastic limit stress and fracture energy, respectively. In order to demonstrate the ability of DCM to simulate cohesive fracture with this simple two-parameter model, several different fracture analyses will be summarized in the following sections. Included will be examples of quasi-static fracture, dynamic crack propagation, and fragmentation.

\subsection{Numerical results}
Multiple numerical tests are carried out to check the efficiency and robustness of the established framework. Quasi-Static and dynamic fracture simulations are performed. 

\subsection{Quasi-Static Fracture}
\begin{figure}[t]
\centering 
{\includegraphics[width=0.75\textwidth]{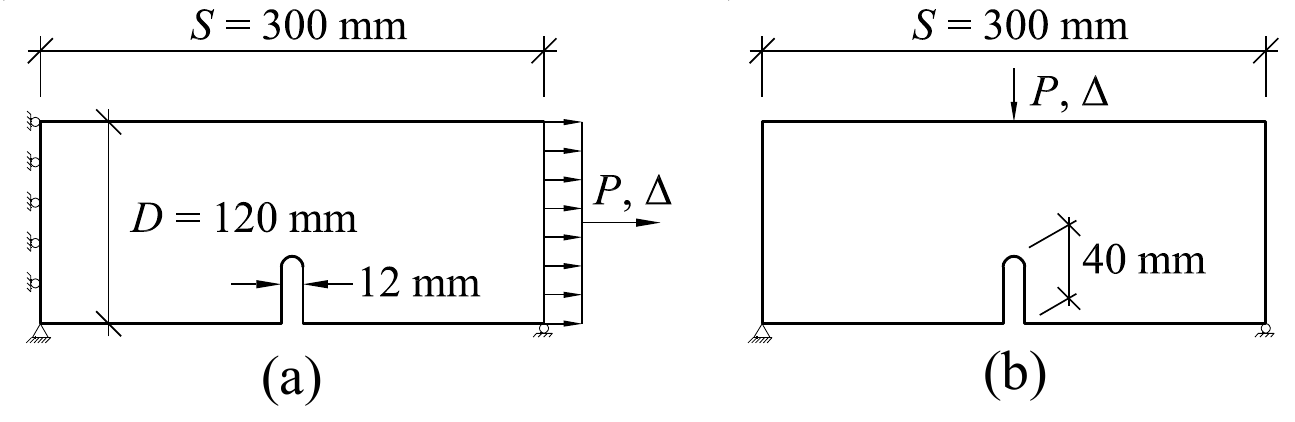}}
\caption{Quasi-Static fracture tests (a) Direct tension test. (b) Three point bending test.}
\label{NotchedSpecimen}
\end{figure}

In this section, the simulation of direct tensile tests and three-point bending tests on notched specimens under quasi-static loading is carried out. The specimens are 120$\times$300 mm rectangular panels with out of plane thickness of 80 mm. The notch is one third of the panel depth, 40 mm, with a width equal to 12 mm. As shown in Figure \ref{NotchedSpecimen}, the notch tip is assumed to be semicircular in order to avoid unrealistic singularities in the stress distribution which might lead to premature crack initiation and propagation. For both the tensile and the three-point bending tests, the assumed material model parameters are: $\varepsilon_t = \sigma_t/E = 6.7 \times 10^{-5}$, tensile strength $\sigma_t = 2$ MPa, Young's modulus $E = 30$ GPa; and material characteristic length $\ell_t = 0.7$ m, which corresponds to a fracture energy of approximately 47 J/m$^2$.

\begin{figure}[h!]
\centering 
{\includegraphics[width=0.8\textwidth]{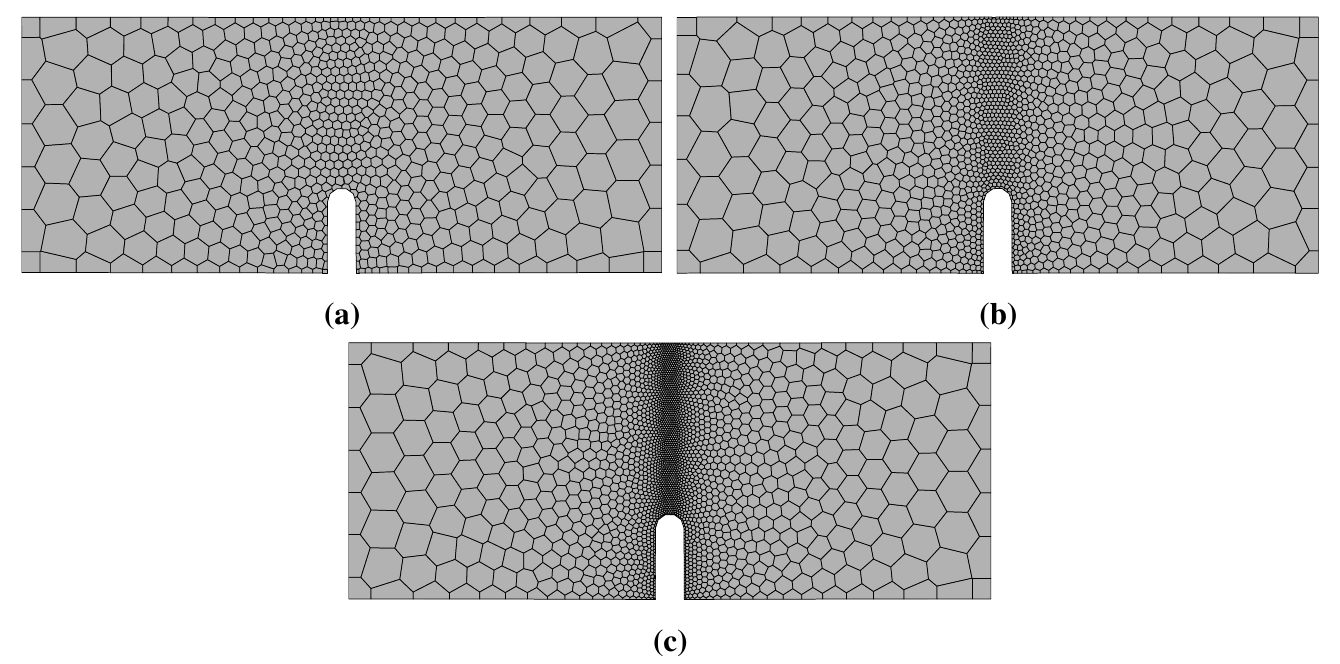}}
\caption{Notched specimens with three different element sizes in the notch zone (a) Coarse (b) Medium (c) Fine.}
\label{NotchedGeoms}
\end{figure}

The direct tensile test is performed by constraining the boundary nodes on the left side of the specimen and by applying an increasing displacement to the nodes on the right side of the specimen, see Figure \ref{NotchedSpecimen}(a). In order to investigate the mesh size dependency of the solution, three different average element size in the zone of the notch of 4, 2, and 1 mm for the coarse, medium, and fine meshes, respectively, are considered and shown in Figure \ref{NotchedGeoms}. 

\begin{figure}[h!]
\centering 
{\includegraphics[width=0.85\textwidth]{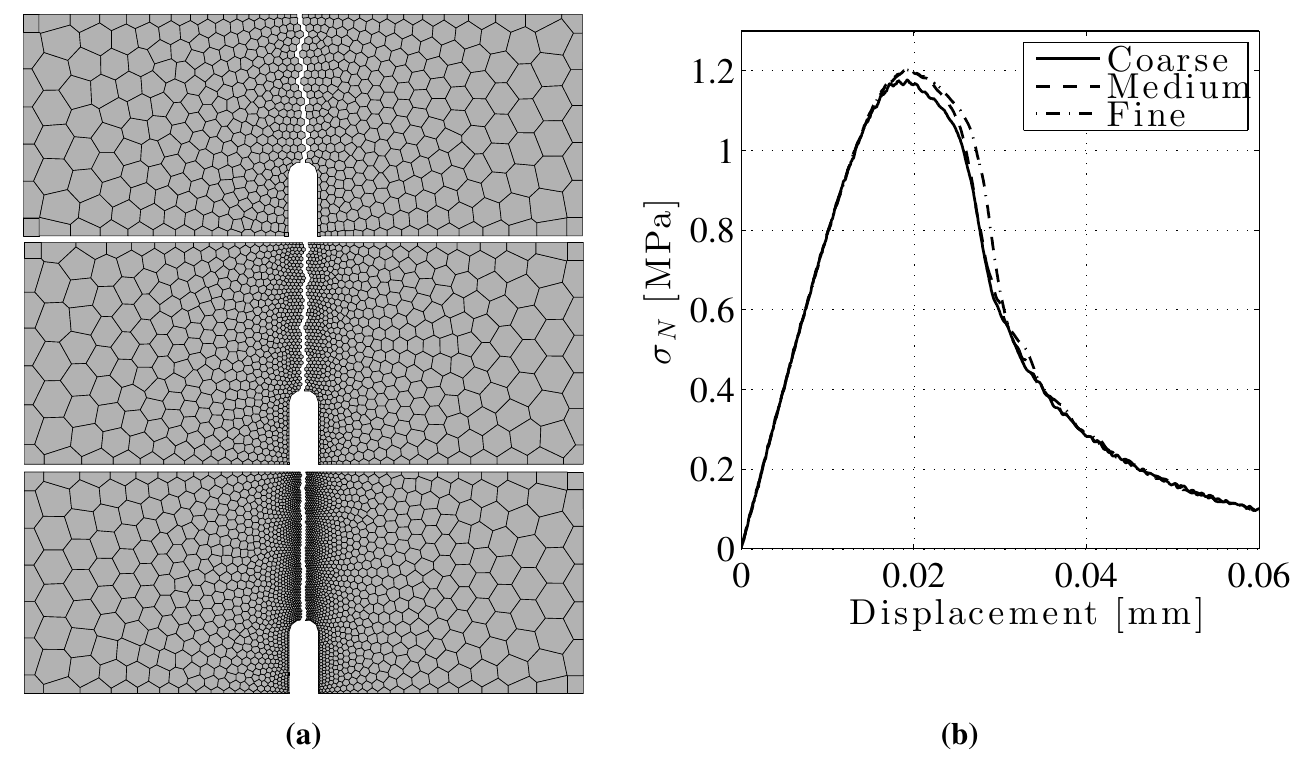}}
\caption{(a) Cracked specimens with different mesh resolutions.(b) Nominal stress versus the right edge displacement.}
\label{NotchedResults-DT}
\end{figure}

Figure \ref{NotchedResults-DT}b reports the nominal stress versus applied displacement curves for three different mesh resolutions. The nominal stress is defined as $\sigma_N = P/Db$, where $P$ is the overall applied load corresponding to a certain displacement. $D$ and $b$ are the specimen width and thickness, respectively. Cracked specimens with different mesh resolutions are plotted in Figure \ref{NotchedResults-DT}a. The crack pattern shown in this figure is relevant to the applied displacement at the end of the loading process (0.06 mm, see Figure \ref{NotchedResults-DT}b). The nodal displacements in Figure \ref{NotchedResults-DT}a have been amplified by 50 to clearly visualize the crack that develops from the notch tip and towards the upper edge of the specimen. All curves reported in Figure \ref{NotchedResults-DT}b have an initial elastic tangent followed by a nonlinear hardening branch up to a certain peak load. Afterwards, the behavior is characterized by a softening response with decreasing load carrying capacity for increasing end displacement. One can observe that the DCM responses for different mesh configurations matches very well which confirms convergence and mesh insensitivity of the DCM framework obtained with the simple regularization in Equation \ref{effmaxstr}.

Spurious mesh sensitivity which is characterized by mesh dependent energy dissipation and lack of convergence is the typical problem of the tensorial constitutive equations for softening materials when not properly regularized \cite{Bazant-CrackBandModel}. Classical numerical techniques such as finite element method that employs these type of constitutive equations result in more brittle response as a finer mesh is utilized. DCM solves this problem by introducing a length type variable $\ell_f$, the local distance between the two nodes, into the vectorial constitutive equations defined on each facet, Equations \ref{effmaxstr}. Therefore, the constitutive equations for any generic facet vary with the edge length between the two corresponding nodes, which yields to a mesh independent global response for different mesh configurations. Fracture energy, energy consumed per unit area of the formed crack, is equal to the area under the curves in Figure \ref{NotchedResults-DT}(b), which means all the three simulations show the same fracture energy. The simulations were not carried out until the applied load dropped to zero, and so only partial estimates of the total energy consumption relevant to complete material separation can be made. For the performed DCM simulations, consumed energies are calculated as 45.7, 46.5, and 46.9 J/m$^2$ for the coarse, medium, and fine meshes, respectively. As noted above the material property used for all simulations was relevant to a fracture energy of $G_t$ $\approx$ 47 J/m$^2$. These positive results provide some quantitative evidence regarding convergence, and proper energy consumption independent of mesh refinement of DCM.

The configuration for the three-point bending test is shown in Figure \ref{NotchedSpecimen}(b). The bottom left corner of the specimen is constrained both vertically and horizontally, whereas the bottom right corner is constrained only in the vertical direction (classical pin-roller boundary conditions for simply supported beams). The load is applied by gradually increasing the displacement of a few boundary nodes located on the top side of the specimen, centered at the specimen axis of symmetry. Similar to the case of direct tension, the obtained results are reported in terms of nominal stress versus applied displacement for the three mesh resolutions. 

\begin{figure}[h!]
\centering 
{\includegraphics[width=0.85\textwidth]{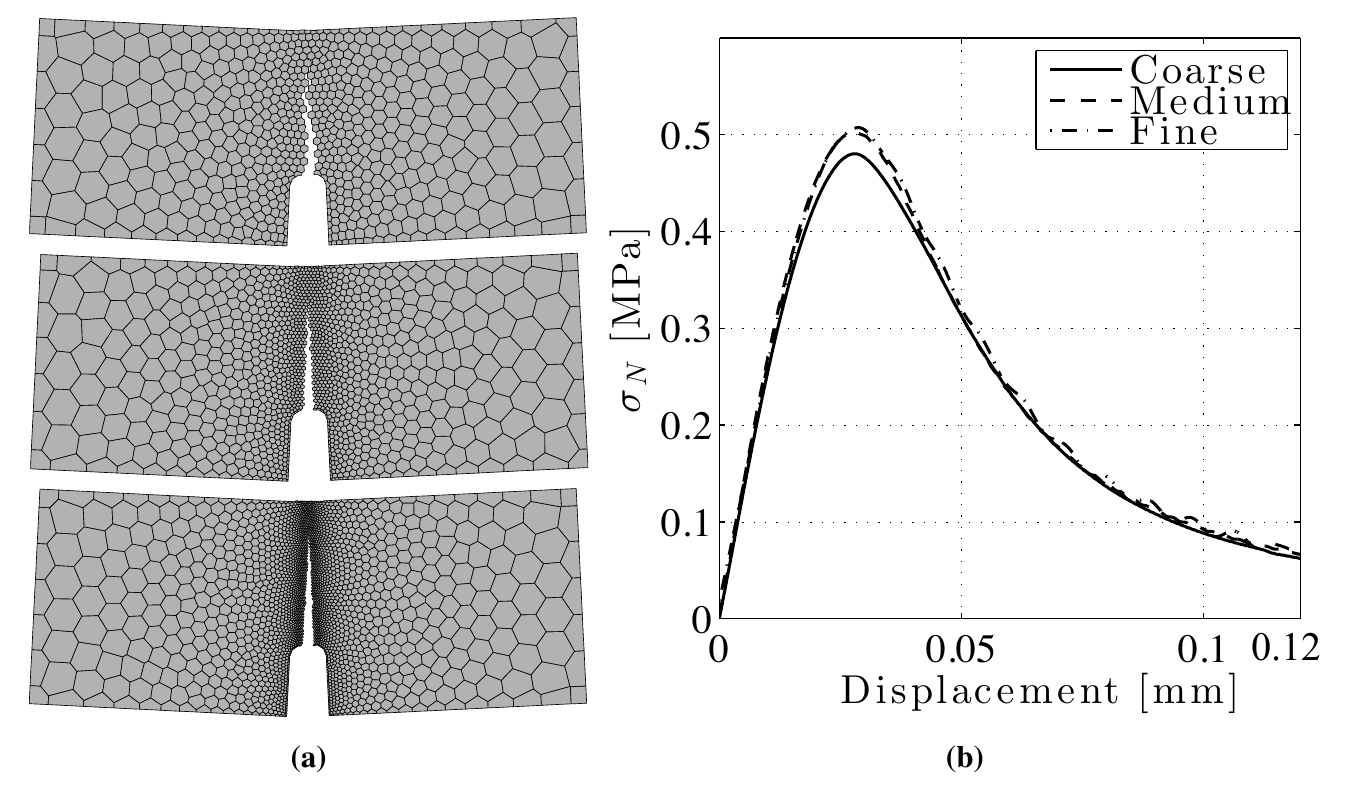}}
\caption{(a) Cracked specimens with different mesh resolutions. (b) Nominal stress versus applied displacement.}
\label{NotchedResults-3pbt}
\end{figure}


Cracked specimens at the end of the loading process (0.12 mm, see Figure \ref{NotchedResults-3pbt}(b)) are depicted in Figure \ref{NotchedResults-3pbt}(a). Due to the bending character of these simulations, the stress profile and the crack opening displacement along the crack are not uniform. Nominal stress versus the applied displacement for the three simulations are plotted in Figure \ref{NotchedResults-3pbt}(b). The results display once again that the DCM solution is convergent upon mesh refinement and spurious mesh sensitivity does not take place.

\subsection{Dynamic Fracture and Fragmentation}
\subsubsection{Bar End-Spalling Fragmentation Test} \label{end-spalling}
\begin{figure}[t]
\centering 
{\includegraphics[width=0.6\textwidth]{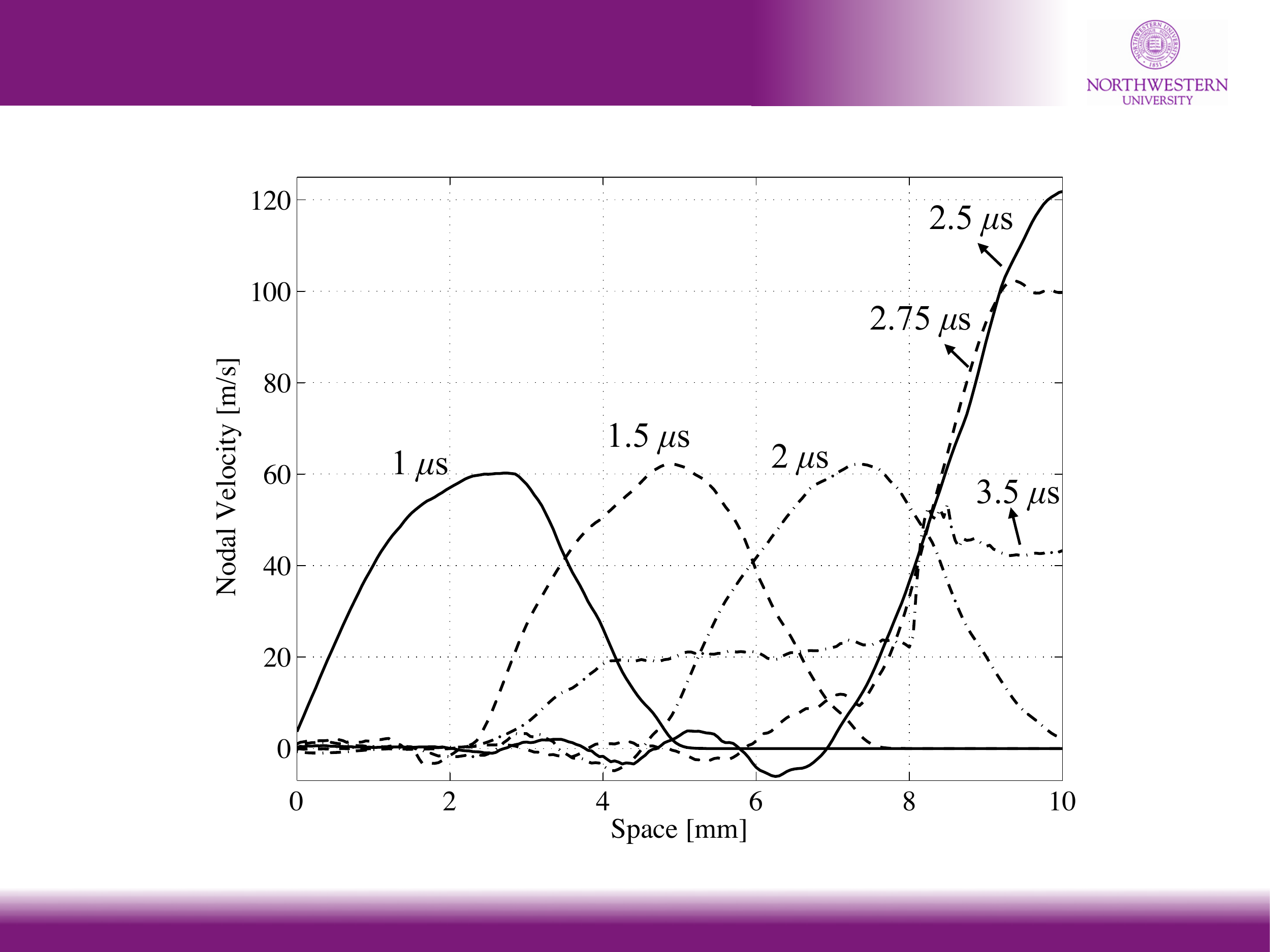}}
\caption{Nodal velocity versus bar longitudinal space at different time steps.}
\label{Bar-Sinus-velTime}
\end{figure}

\begin{figure}[h!]
\centering 
{\includegraphics[width=0.8\textwidth]{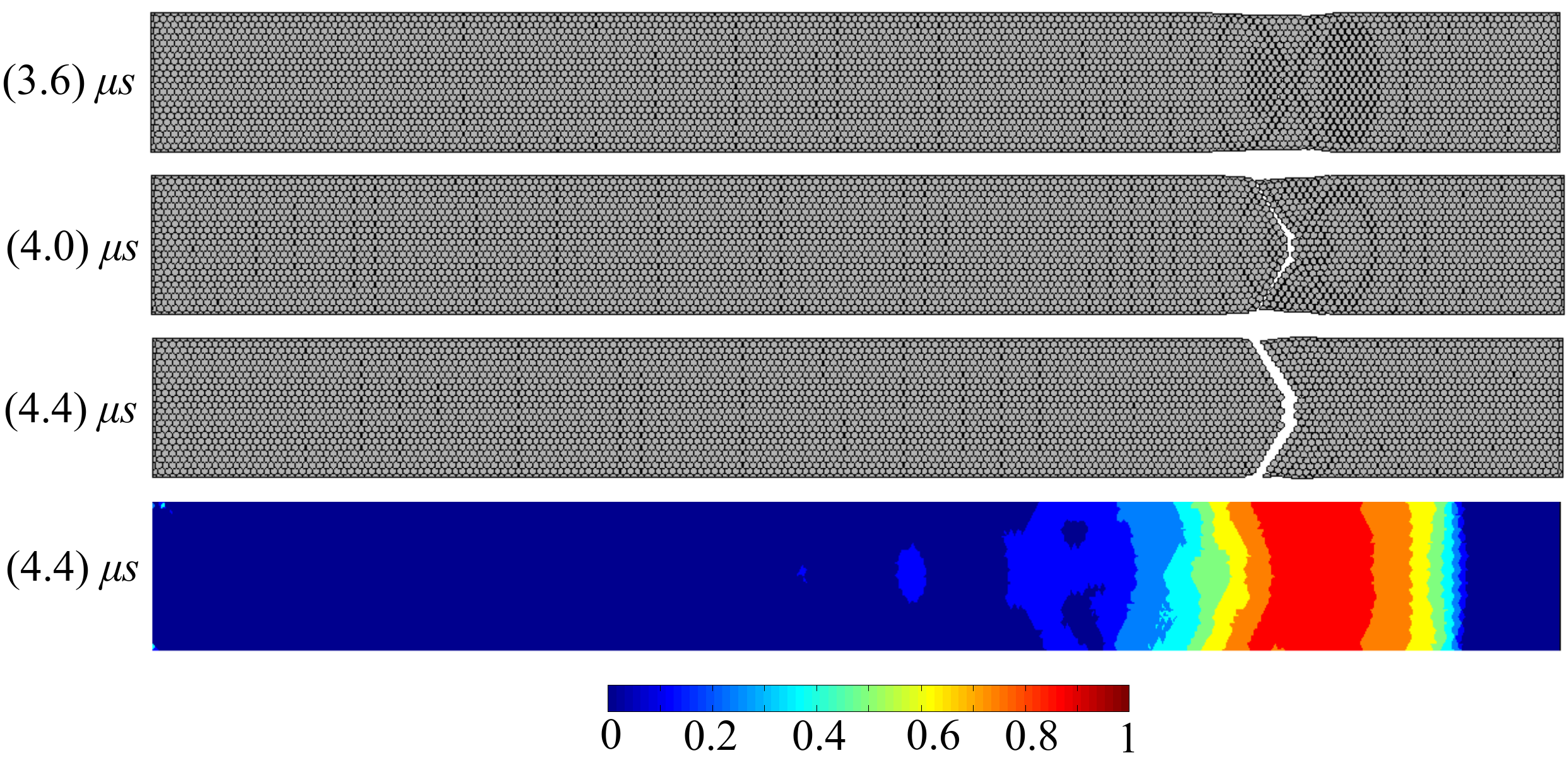}}
\caption{Fracture pattern and damage variable $D_f$ contour of a bar under sinusoidal velocity impulse at failure.}
\label{crack-bar-impulse}
\end{figure}

An end-spalling fragmentation experiment is carried out on a two dimensional bar of 1 mm height $\times$ 10 mm length and 1 mm thickness under plane stress condition. The bar is subjected to a sinusoidal velocity impulse with peak value of $60$ m/s and  1 $\mu$s duration time. In this study, the assumed material properties are: $\sigma_t = 844$ MPa, $E = 190$ GPa, material characteristic length $\ell_t = 0.012$ m which corresponds to a fracture energy of 2.2 $\times$ $10^4$ J/m$^2$, material density $\rho = 8000$ $\text{kg/m}^\text{3}$, and the Poisson's ratio $\nu = 0.30$. The horizontal velocity component of the nodes located vertically at mid-height of the bar is plotted along the bar (Figure \ref{Bar-Sinus-velTime}). One can see that the velocity impulse travels undisturbed across the bar, see time steps 1, 1.5, and 2 $\mu$s, while the magnitude of the wave doubles as it reaches the right end due to the interaction with the free-end boundary condition, see time step 2.5 $\mu$s in Figure \ref{Bar-Sinus-velTime}. One can calculate that the velocity impulse moves across the specimen at an approximate velocity of 5000 m/s, which corresponds closely to the basic equation for one dimensional wave propagation velocity $v = \sqrt{E/\rho}$ which results in a velocity of 4874 m/s. In terms of stress, the traveling wave applies a compressive stress on the bar before arriving at the free-end, where the stress reverses from compression to tension. The generated tensile stress overcomes the tensile strength of the material, and end-spalling fragmentation begins. One can see that the wave moving back through the bar is no more sinusoidal which is due to the engendered material nonlinearity, see time steps 2.75 and 3.5 $\mu$s. More detailed analysis of the problem under consideration reveals that actually a bi-axial strain state is generated through the specimen due to the Poisson's effect results in the presence of lateral straining. In turn, this leads to inclined principal strain directions and, consequently, inclined cracks (see Figure \ref{crack-bar-impulse}) since crack initiation and propagation are simulated with the strain dependent damage model discussed earlier.

\begin{figure}[h!]
\centering 
{\includegraphics[width=0.8\textwidth]{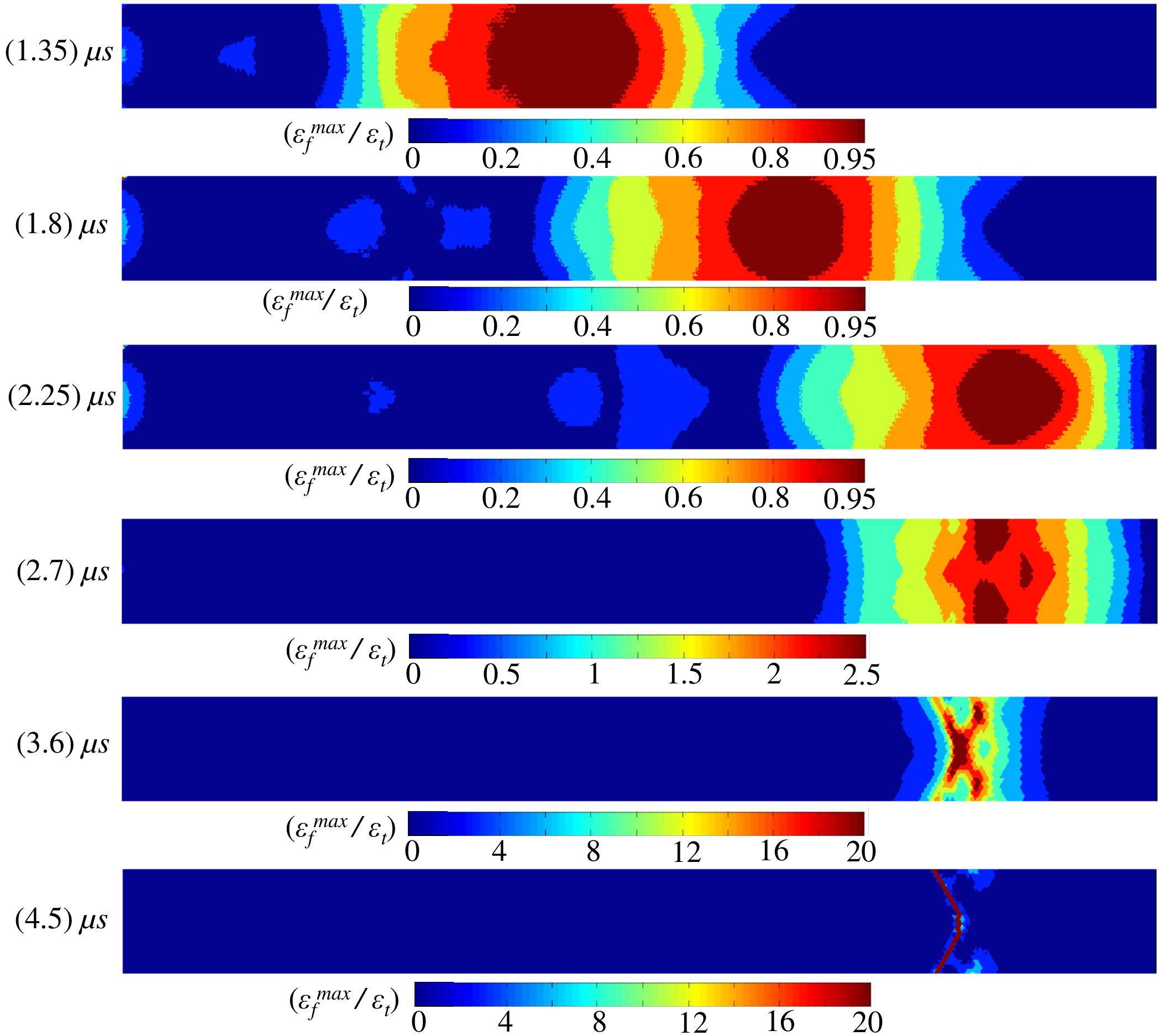}}
\caption{$\varepsilon_{f}^{max}/\varepsilon_{t}$ ratio contour at different time instants.}
\label{bar-eps-ratio}
\end{figure}

Figure \ref{crack-bar-impulse} shows the fracture pattern of the bar at different time steps namely 3.6, 4, and 4.4 $\mu$s. One can see that a localized fracture takes place at the bar end and evolves into total separation of the right end once the fracture energy of the material is completely overcome. Contour of the damage parameter $D_f$ is also illustrated in Figure \ref{crack-edge-impulse}, which confirms the fracture pattern occurred at the bar end. In addition, one can notice that the crack propagates vertically at the center of the bar, while it deviates as it moves towards the cross section edges. This can be explained by the fact that the bi-axial effect is more pronounced over the areas away from the center of the bar. In Figure \ref{bar-eps-ratio}, maximum effective strain experienced by each facet $\varepsilon_{f}^{max}$ normalized by $\varepsilon_t = \sigma_t/E = 4.4 \times 10^{-3}$ is plotted at different time instants. 1.35, 1.8, and 2.25 $\mu s$ are instants during which the compressive wave travels through the bar before reaching the free-end, while 2.7, 3.6, and 4.5 $\mu s$ are after the signal reaches the free-end and leads to a tensile wave. One can see that at 1.35, 1.8, and 2.25 $\mu s$, $\varepsilon_{f}^{max}/\varepsilon_t < 1$ for all facets, which implies that all facets stay in the elastic regime. At 2.7 $\mu s$, which is just after the signal reaches the free-end, and the compressive wave converts into tensile one, $\varepsilon_{f}^{max}/\varepsilon_t > 1$ at the end of the bar where nonlinearity starts to develop. At 3.6 and 4.5 $\mu s$, damage localizes, and $\varepsilon_{f}^{max}/\varepsilon_t$ contour corresponds to the specimen fracture pattern depicted in Figure \ref{crack-bar-impulse}.

\begin{figure}[t]
\centering 
{\includegraphics[width=0.85\textwidth]{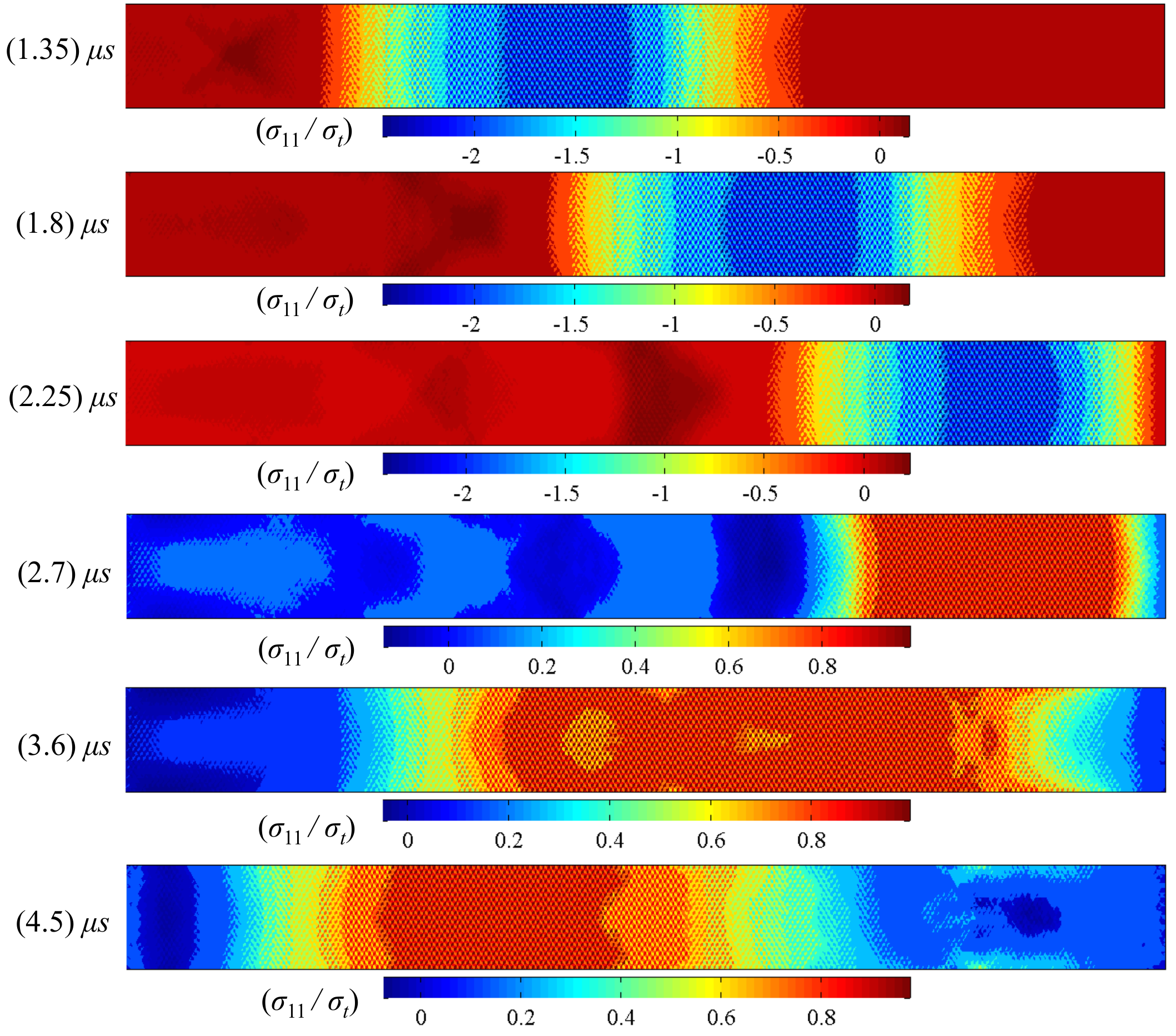}}
\caption{$\sigma_{11}/\sigma_{t}$ ratio contour at different time instants.}
\label{bar-sig-ratio}
\end{figure}

Ratio of the horizontal component facet stress tensor $\sigma_{f,11}$ to $\sigma_t$ is plotted in Figure \ref{bar-sig-ratio} at the same time instants as the ones considered in Figure \ref{bar-eps-ratio}. One can clearly see the propagation of the compressive wave through the specimen at 1.35, 1.8, and 2.25 $\mu s$, and its conversion to tensile wave at 2.7 $\mu s$. At 4.5 $\mu s$, it can be observed that the stress value on the facets around which fracture takes place is approximately zero, which corresponds to the bar splitting type of failure pattern. 

\subsubsection{Edge Cracked Plate under Velocity Impulse}
\begin{figure}[h!]
\centering 
{\includegraphics[width=0.6\textwidth]{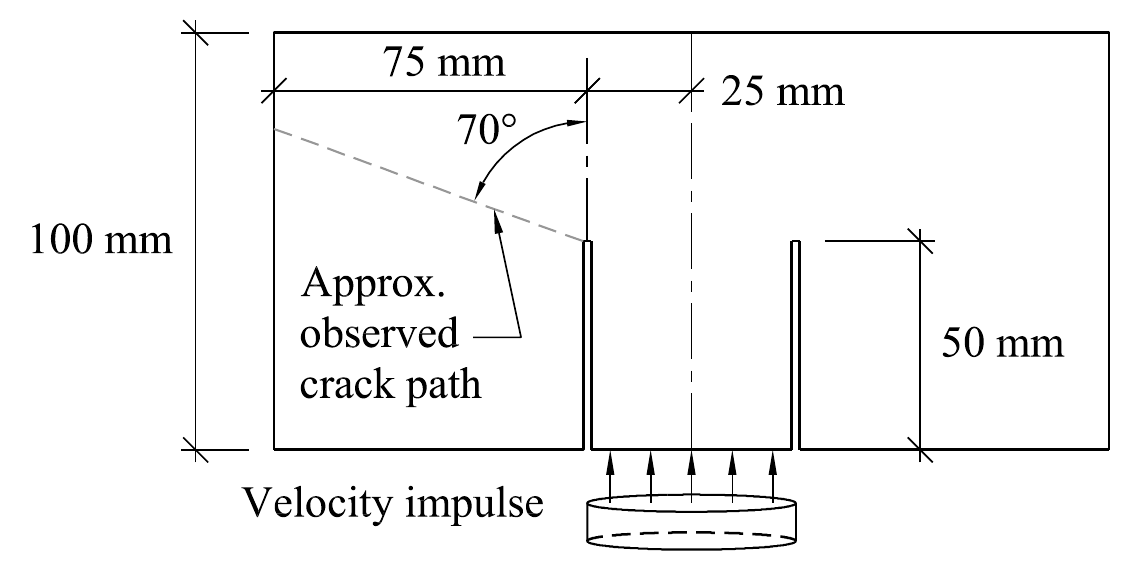}}
\caption{Experimental set up for edge cracked plate under velocity impulse.}
\label{crack-edge-impulse}
\end{figure}

In this section, a classical dynamic crack propagation test is simulated. The reference experimental data is relevant to maraging steel \cite{Kalthoff-1}, which shows high tensile strength and brittle behavior when subjected to high strain rate. A schematic representation of the test configuration is shown in Figure \ref{crack-edge-impulse}, in which one can see a projectile impacting the central part of an unrestrained double notched specimen. The plate has a 10 mm out-of-plane thickness. Plane stress condition can be assumed for the DCM analysis. By using the symmetry of the problem, half of the specimen is modeled and appropriate boundary conditions, horizontal nodal displacement and nodal rotation equal to zero, are enforced on the line of symmetry. Kalthoff and Wrinkler \cite{Kalthoff-1} investigated the effect of the projectile velocity on the failure mechanism: a brittle failure with a crack at an angle of $-70^{\circ}$ was observed for the case of low impact velocity (32 m/s), see Figure \ref{crack-edge-impulse}. In the current example, a velocity of 16 m/s is applied at the impacted nodes, and this impulse is kept constant to the end of the simulation. The velocity of 16 m/s is selected because the elastic impedance of the projectile and the specimen are considered to be equal. Material properties considered in the DCM simulations are the same as the ones used in Section \ref{end-spalling}.

\begin{figure}[h!]
\centering 
{\includegraphics[width=0.95\textwidth]{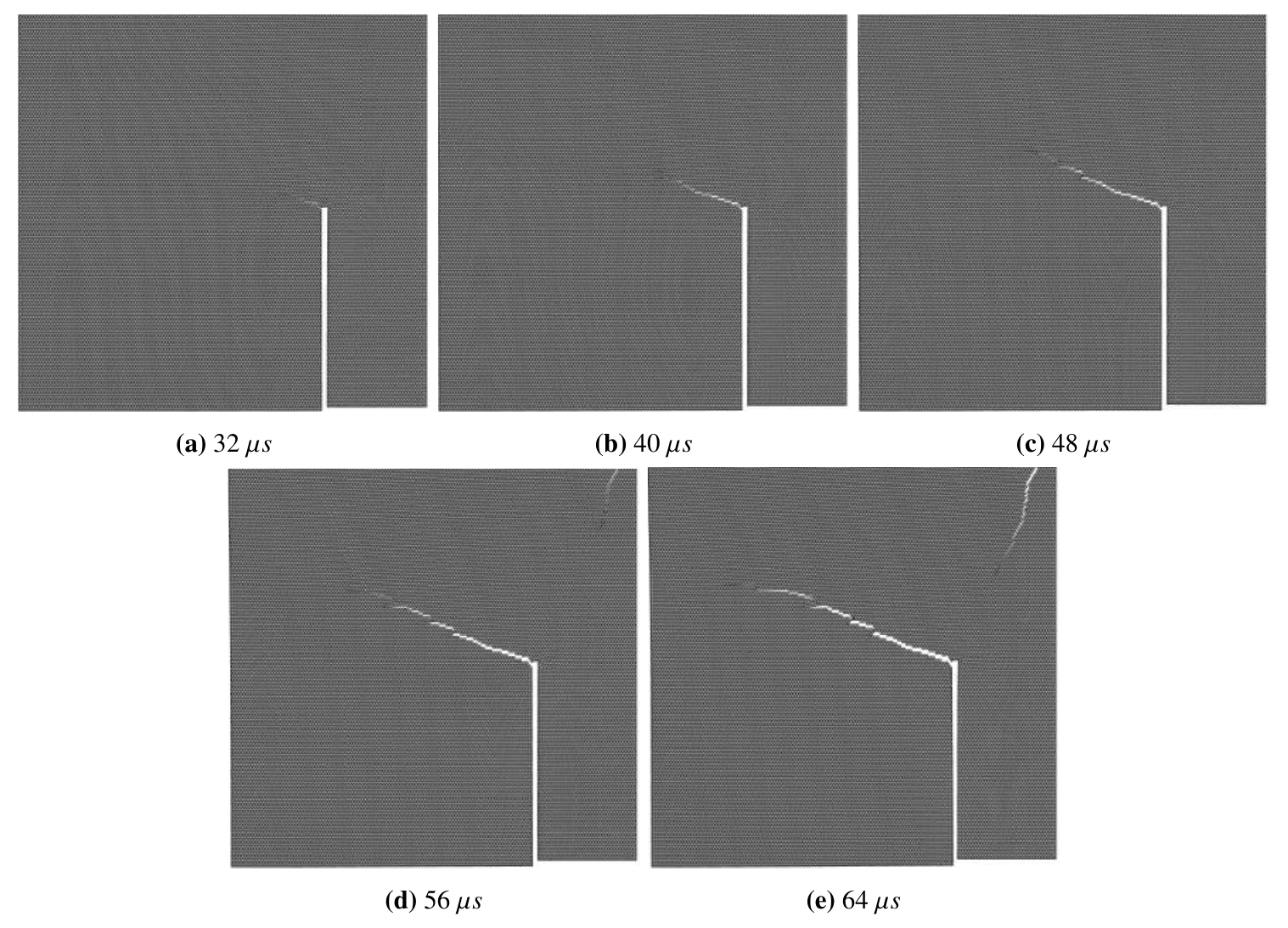}}
\caption{Fracture mode of the fine mesh edge cracked plate under impulsive load test.}
\label{ECPF-Frac}
\end{figure}

Belytschko et al. \cite{Belytschko-1} simulated this experiment using the continuum based model XFEM for quasi-brittle fracture, which is considered here as reference to discuss the DCM performance. To investigate the mesh dependency of the DCM results, a fine and a coarse mesh with element edge of $\sim$0.65 mm (50573 elements) and $\sim$1 mm (22437 elements) are considered. The initial vertical notch is simulated with 1.5 mm width, and the time step used in the explicit integration scheme is 0.02 $\mu$s. The  velocity impulse applied to the DCM boundary nodes generates a compressive wave in the central part of the specimen, which propagates until it reaches the notch tip. At this point, significant shear strains develop leading to high principal tensile strains and crack initiation at the left side of the notch tip. Subsequently, crack propagates towards the left boundary of the specimen.

\begin{figure}[t]
\centering 
{\includegraphics[width=0.9\textwidth]{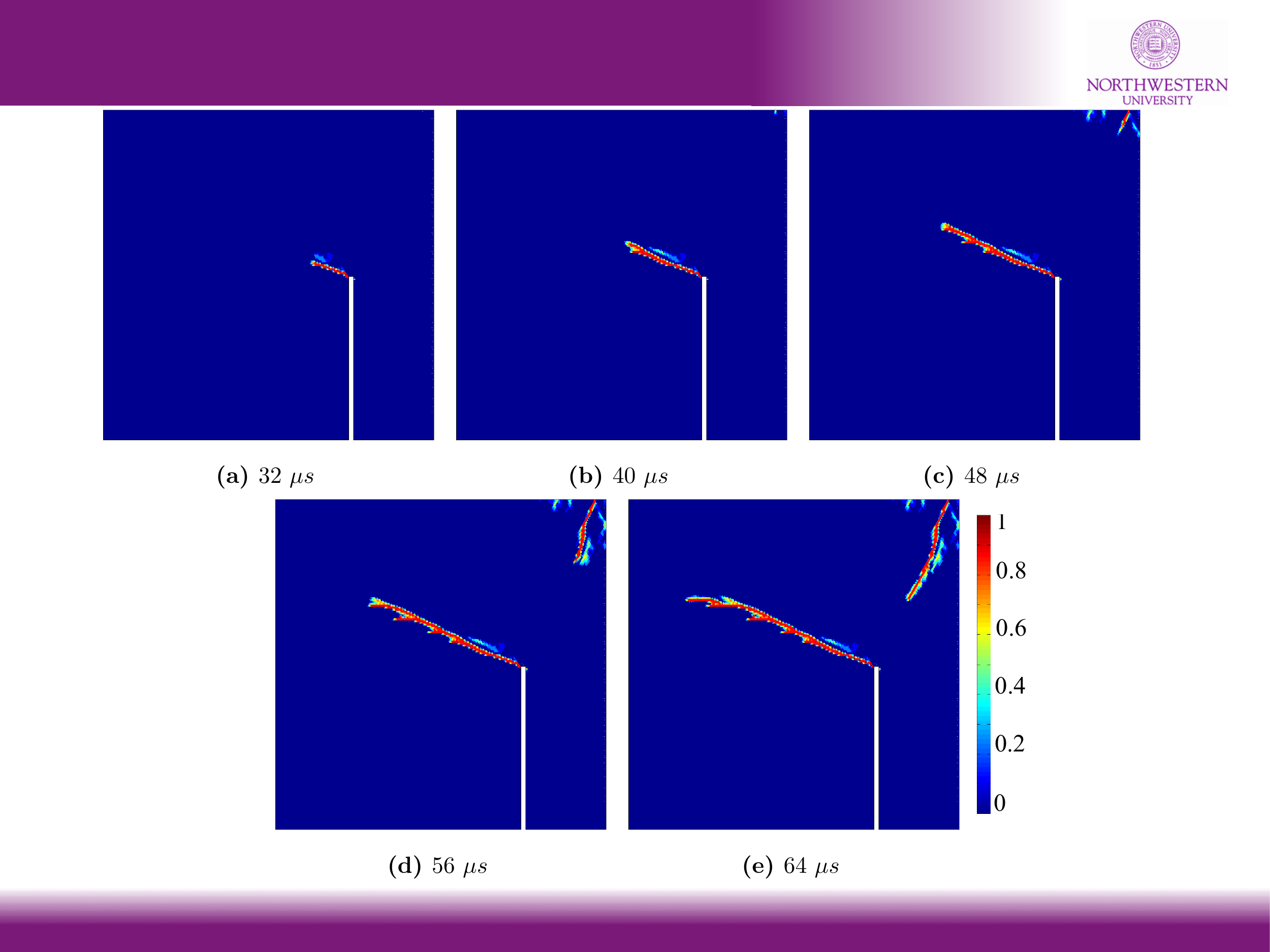}}
\caption{Damage parameter contour of the fine mesh edge cracked plate under impulsive load test.}
\label{ECPF-Damage}
\end{figure}

Figure \ref{ECPF-Frac} shows crack initiation and propagation in different time steps for the fine mesh case. The average crack propagation angle with the horizontal axis at the time step 56 $\mu$s is $69^\circ$ which compares very well with experimental result $70^\circ$. At this time, a localized damage which leads to fracture takes place on the top right boundary of the specimen, and the generated crack propagates towards the notch tip. This is due to the reflection of the compressive wave from the top right boundary and is also reported by  Belytschko et al. \cite{Belytschko-1} in their XFEM simulations. The propagating crack tends to become horizontal at the end of the simulation, see Figure \ref{ECPF-Frac}e, as the localized fracture occurs and propagates on the top right boundary. This can also be related to the strain based failure criteria employed in DCM model and the simple damage model used in constitutive behavior of the material. Damage coefficient $D_f$ contours of the fine mesh simulation are plotted in Figure \ref{ECPF-Damage}, which clearly shows two highly localized damaged areas corresponding to the fracture pattern depicted in Figure \ref{ECPF-Frac}.  

\begin{figure}[t]
\centering 
{\includegraphics[width=0.95\textwidth]{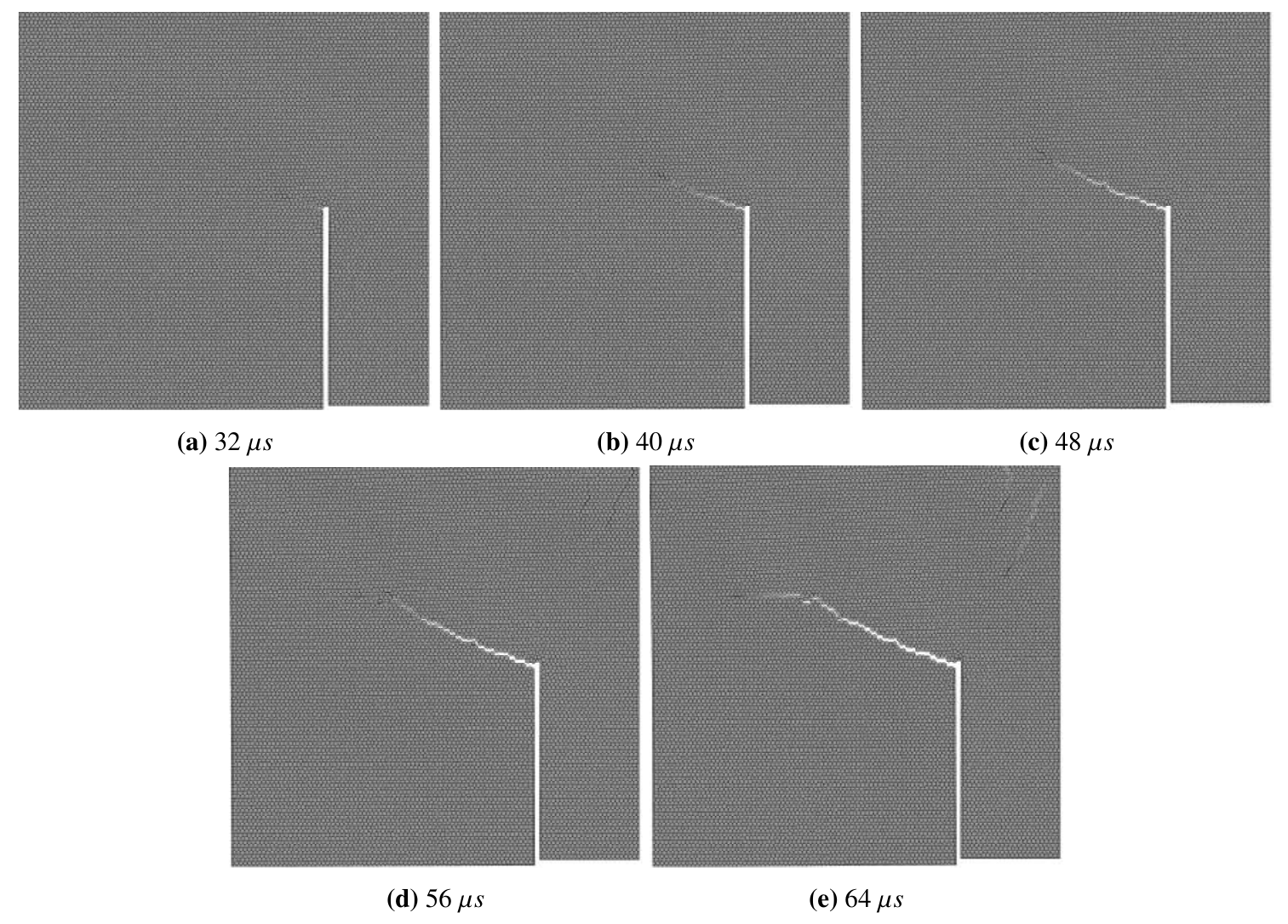}}
\caption{Fracture mode of the coarse mesh edge cracked plate under impulsive load test.}
\label{ECPC-Frac}
\end{figure}

\begin{figure}[t]
\centering 
{\includegraphics[width=0.9\textwidth]{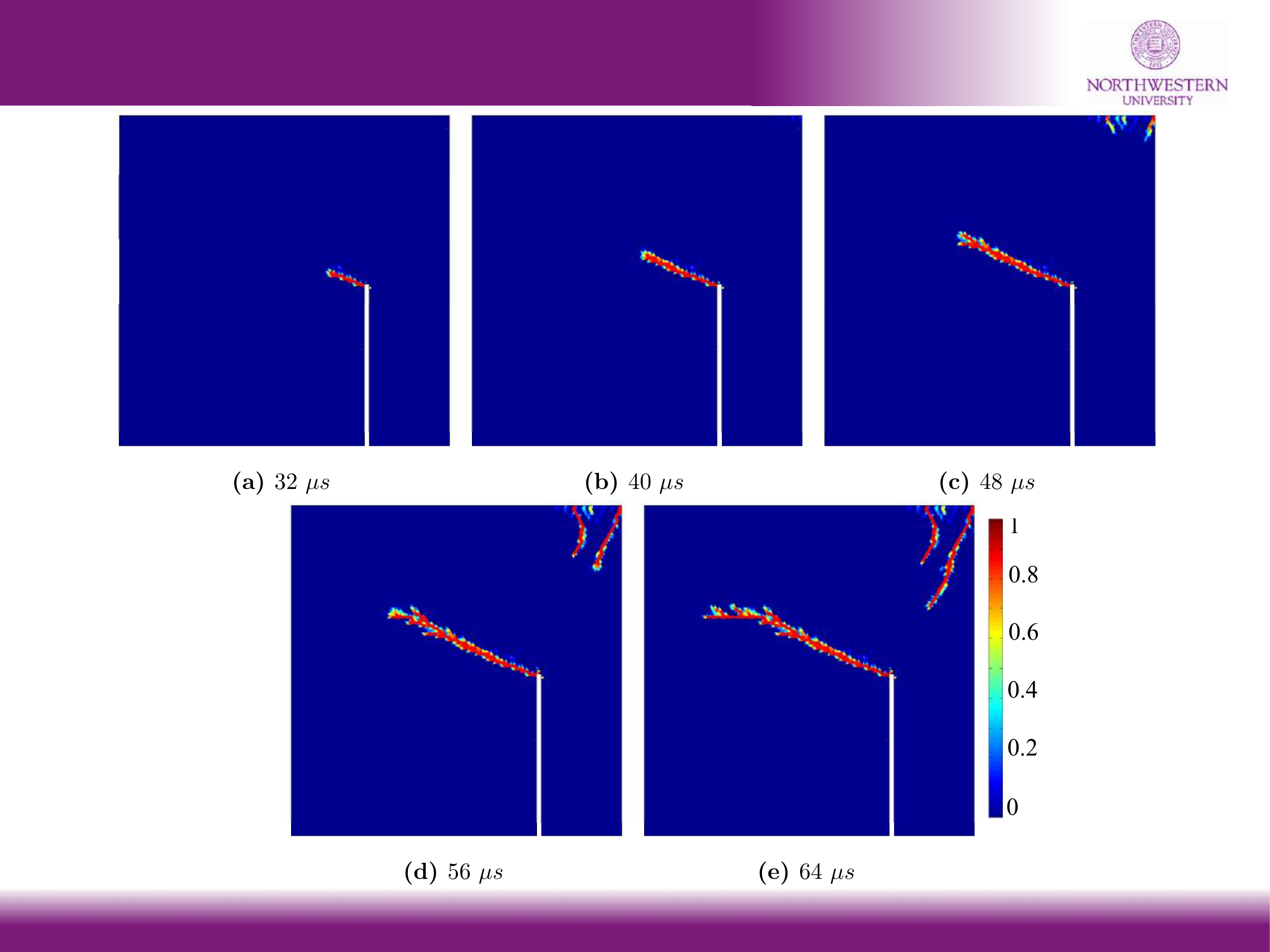}}
\caption{Damage parameter contour of the coarse mesh edge cracked plate under impulsive load test.}
\label{ECPC-Damage}
\end{figure}

Figures \ref{ECPC-Frac} and \ref{ECPC-Damage} shows the fracture pattern and damage coefficient contour at different time steps for the coarse mesh simulations, respectively, which agrees well with the fine mesh results. Average crack propagation angle with the horizontal axis at the time step 56 $\mu$s is $68^\circ$, and the crack tends to propagate horizontally as fracture occurs and develops from the top right boundary. DCM performs more accurately compared to the XFEM \cite{Belytschko-1} in terms of predicting the crack propagation angle. However, the crack does not develop on the same path to the end of the test, which is captured by the XFEM simulations. In addition, DCM is able to capture micro cracks developing from the main crack faces, as it can be seen in Figures \ref{ECPF-Damage} and \ref{ECPC-Damage}, while this is not captured by other approaches. It is also worth noting that the computational cost of methods like XFEM increases as the crack propagates because additional DOFs must be inserted to capture the displacement discontinuity. This is not the case for DCM which is characterized by the same number of DOFs in the elastic and fracturing regimes. 

\subsubsection{Dynamic Crack Branching}
\begin{figure}[h!]
\centering 
{\includegraphics[width=0.35\textwidth]{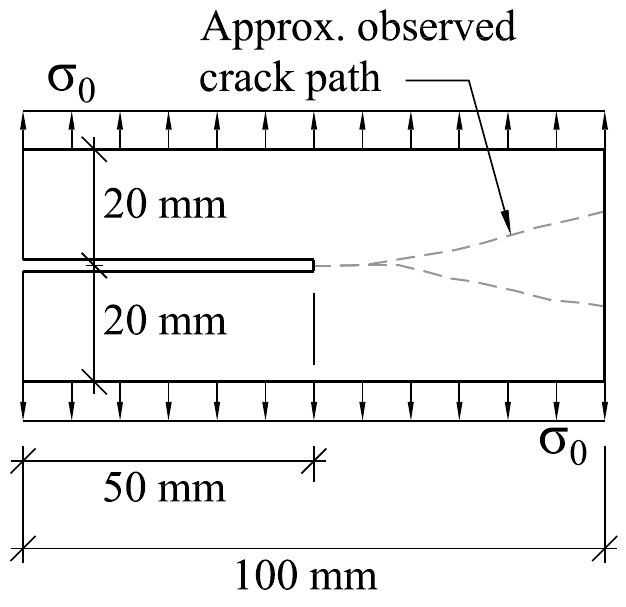}}
\caption{Dynamic crack branching: test geometry and experimental results.}
\label{Branch-Geom}
\end{figure}

A final benchmark fracture problem simulated by DCM in this section involves dynamic crack propagation and crack branching. Figure \ref{Branch-Geom} shows a schematic representation of the test configuration. A pre-notched rectangular panel is subjected to a uniform traction applied as a step function on the two edges parallel to the notch. This experiment has been simulated computationally by other authors \cite{Song-1, Xu-1}, and related experimental results were reported by different researchers \cite{Ramulu-1,Ravi-Chandar-1}. Ramulu and Kobayashi \cite{Ramulu-1} observed experimentally that a major crack starts to propagate from the notch tip to the right, which branches into two cracks at a certain point during the experiment, see Figure \ref{Branch-Geom} for the sketch of the experimental result. The DCM parameters used in this test are: $\sigma_t = 3.1$ MPa, $E = 32$ GPa; material characteristic length $\ell_t = 0.02$ m, which corresponds to a fracture energy of approximately 3 J/m$^2$; material density $\rho = 2500$ $\text{kg/m}^\text{3}$; Poisson's ratio $\nu = 0.2$. The applied traction is $\sigma_0 = 1$ MPa. 

\begin{figure}[t]
\centering 
{\includegraphics[width=0.95\textwidth]{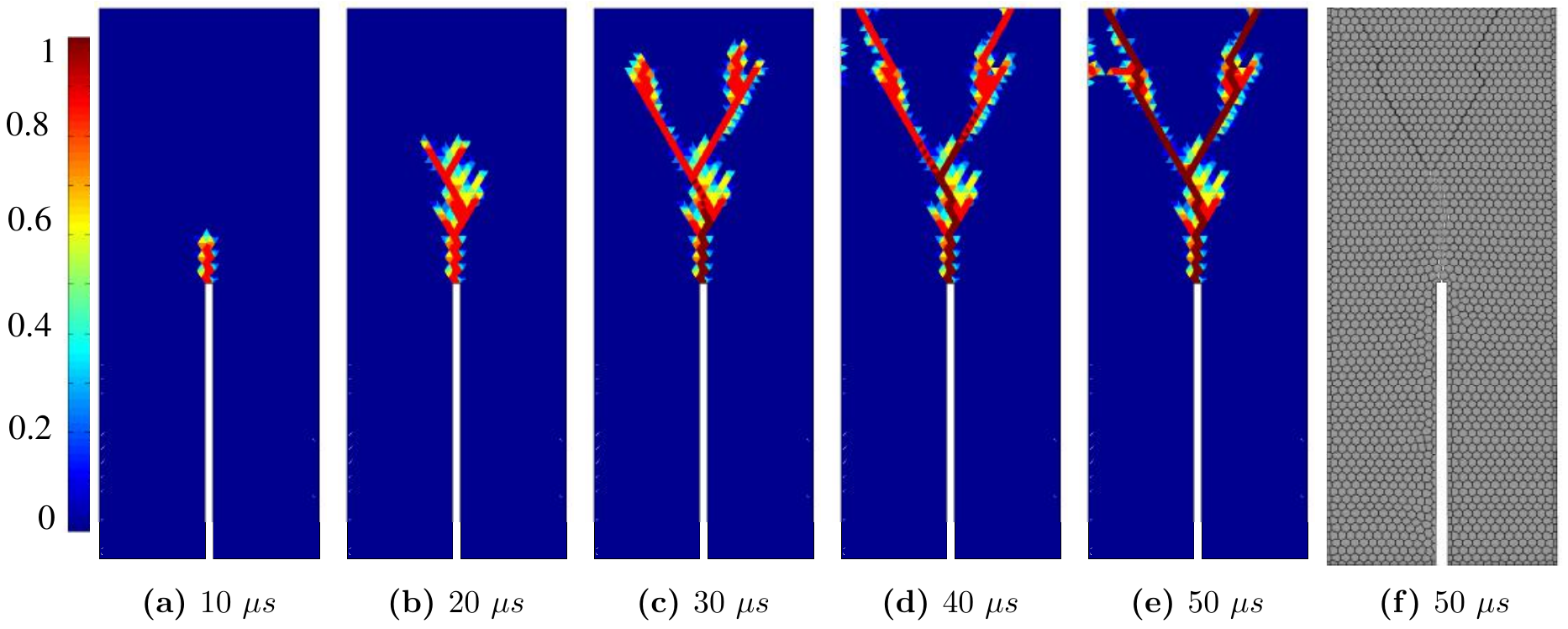}}
\caption{(a-e) Damage of the crack branching test at different time steps. (f) Fracture mode of the specimen at the end of the simulation.}
\label{Branch-Damage}
\end{figure}

Crack initiation and propagation at different time steps of the DCM simulation is illustrated in Figure \ref{Branch-Damage}(a-e) through the damage parameter contours. One can see that the crack starts to propagate from the notch tip parallel to the symmetry axis of the configuration on a straight path for a short distance, and it branches into two cracks subsequently. The deformed configuration of the DCM simulation is plotted in Figure \ref{Branch-Damage}f, which agrees well with the experimental results. Experimental observations also report that before the main branching occurs, minor branches emerge from the main crack but only propagate on a short distance \cite{Ramulu-1}. DCM is able to capture these minor branches which can be seen in the damage variable contours in Figure \ref{Branch-Damage}.

\section{DCM Versus Particle Methods}
DCM and classical particle models are basically governed by the same set of algebraic equations expressing compatibility and equilibrium. This naturally follows from the adoption of rigid body kinematics which is common to the two approaches. The difference between the two methods lies in the formulation of the constitutive law, namely in the relationship between facet tractions $\bold{t}_f$ and facet openings $\bold{w}_f$. For classical particle models, in which rigid particles are connected by springs, this relationship can be expressed as $t_{fN} = E_N w_{fN}/\ell_f$ and $t_{fM} = E_T w_{fM}/\ell_f$, where $E_N$ and $E_T$ are the normal and tangential elastic stiffnesses, respectively. DCM constitutive laws, Equation \ref{NTConstDamag}, are consistent to the ones for particle models if one sets $E_V = E_D$ and $E_T = E_D$. These conditions correspond to an elastic material with zero Poisson's ratio \cite{Bolander-1}. By properly setting the ratio between the normal and tangential stiffnesses, particle models can simulate an “average” non-zero Poisson’s
￼￼￼
ratio (average in the sense that Poisson's ratio is defined by analyzing a finite, as opposed to an infinitesimal, volume of material). In this case, however, particle models feature an intrinsic heterogeneous response even for load configurations that produce uniform strain fields according to continuum theory \cite{Bolander-1}. In conclusion, for non-zero Poisson's ratio the two formulations are fundamentally different and the key difference is that DCM accounts for the orthogonality of the deviatoric and volumetric deformation modes while classical particle models do not.

It must be mentioned here that the heterogeneous response of particle models is not necessarily a negative property and, actually, it is critical for their ability to handle automatically strain localization and crack initiation \cite{cusatis-ldpm-1,cusatis-ldpm-2}. It must be kept in mind, however, that in this case the size of the discretization cannot be user-defined but must be linked to the actual size of the material heterogeneity. Only under this condition can one consider the heterogeneous response of particle models to be a representation of the actual internal behavior of the material rather than a spurious numerical artifact.

\section{Conclusions}
In this paper, the formulation of the Discontinuous Cell Method (DCM) has been outlined. A convergence study in the elastic regime shows that DCM converges to the exact continuum solution with a convergence rate that is comparable to that of constant strain finite elements, but with accuracy that is one order of magnitude higher. In addition, numerical simulations show that DCM, with a simple two parameter isotropic damage model, can simulate cohesive fracture propagation without the drawbacks of standard finite elements, such as spurious mesh sensitivity, and without the complications of most recently formulated computational techniques. In addition, DCM successfully simulated the crack branching which is observed in the experiment of a benchmark problem. Finally, DCM can simulate the transition from localized fracture to fragmentation without mesh entanglement typical of finite element approaches.

\noindent\textbf{ACKNOWLEDGMENTS} \\
This material is based upon work supported by the National Science Foundation under grant no. CMMI-1435923.

\newpage
\begin{appendices}
\numberwithin{equation}{section}

\section{Facet Strain Tensor Calculation} \label{facet strain tensor}
Derivation of facet strain tensor components $\varepsilon_{ij}$ in terms of facet normal $\varepsilon_{fN}$, tangential $\varepsilon_{fM}$ and volumetric strains $\varepsilon_{V}$ are explained in this section for plane strain and stress problems. In both cases, one should solve a system of three algebraic equations with three unknowns.

\subsection{Plane Strain Problem}
In plane strain problems, out of plane strain component $\varepsilon_{33} = 0$, and one should consider following system of equations 
\begin{eqnarray}\label{PEstrainTens-1}
\begin{split}
&\varepsilon_{V} = \frac{1}{3}\varepsilon_{ii} = \frac{1}{3}(\varepsilon_{11} + \varepsilon_{22}) \\
&\varepsilon_{fN} = N_{ij} \varepsilon_{ij} = N_{11}\varepsilon_{11} + 2 N_{12}\varepsilon_{12} + N_{22}\varepsilon_{22} \\
&\varepsilon_{fM} = M_{ij} \varepsilon_{ij} = M_{11}\varepsilon_{11} + 2 M_{12}\varepsilon_{12} + M_{22}\varepsilon_{22}
\end{split}
\end{eqnarray}

where $N_{ij} = n_i n_j$ and $M_{ij} = n_i m_j$ are the two projection tensors that are calculated for each facet using its unit normal $n_i$ and tangential $m_i$ vectors, subscript $f$ is dropped for simplicity. Solution of the above system of equations will yield to the following expressions for strain tensor components:
\begin{eqnarray}\label{PEstrainTens-2}
\begin{split}
&\varepsilon_{11} = \frac{\varepsilon_{M} N_{12} - \varepsilon_{N} M_{12} + 3 \varepsilon_{V}(M_{12}N_{22} - N_{12}M_{22})} {M_{11}N_{12} + M_{12}N_{22} - M_{12}N_{11} - M_{22}N_{12}} \\ 
\\
&\varepsilon_{12} = \frac{\varepsilon_{N} (M_{11} - M_{22}) + \varepsilon_{M} (N_{22} - N_{11}) - 3 \varepsilon_{V} (M_{11}N_{22} - N_{11}M_{22})} {2(M_{11}N_{12} + M_{12}N_{22} - M_{12}N_{11} - M_{22}N_{12})} \\ 
\\
&\varepsilon_{22} = \frac{\varepsilon_{N} M_{12} - \varepsilon_{M} N_{12} + 3 \varepsilon_{V}(M_{11}N_{12} - N_{11}M_{12})} {M_{11}N_{12} + M_{12}N_{22} - M_{12}N_{11} - M_{22}N_{12}}
\end{split}
\end{eqnarray}

These quantities are then used to calculate the strain tensor eigenvalues. 

\subsection{Plane Stress Problem}
for the case of plane stress problems, out of plane strain component $\varepsilon_{33} = -\nu(\varepsilon_{11}+\varepsilon_{22})/(1-\nu)$ should be taken into account. Therefore, the first equation in the system of equations \ref{PEstrainTens-1} should be revised as $\varepsilon_{V} = \varepsilon_{ii}/3 = (\varepsilon_{11} + \varepsilon_{22} - \nu(\varepsilon_{11}+\varepsilon_{22})/(1-\nu))/3$, while the two other equations stay the same.
\begin{eqnarray}\label{PEstressTens-2}
\begin{split}
&\varepsilon_{11} = \frac{\varepsilon_{M} N_{12} - \varepsilon_{N} M_{12} + 3(1-\nu)/(1-2\nu) \varepsilon_{V}(M_{12}N_{22} - N_{12}M_{22})} {M_{11}N_{12} + M_{12}N_{22} - M_{12}N_{11} - M_{22}N_{12}} \\ 
\\
&\varepsilon_{12} = \frac{\varepsilon_{N} (M_{11} - M_{22}) + \varepsilon_{M} (N_{22} - N_{11}) - 3 (1-\nu)/(1-2\nu) \varepsilon_{V} (M_{11}N_{22} - N_{11}M_{22})} {2(M_{11}N_{12} + M_{12}N_{22} - M_{12}N_{11} - M_{22}N_{12})} \\ 
\\
&\varepsilon_{22} = \frac{\varepsilon_{N} M_{12} - \varepsilon_{M} N_{12} + 3(1-\nu)/(1-2\nu) \varepsilon_{V}(M_{11}N_{12} - N_{11}M_{12})} {M_{11}N_{12} + M_{12}N_{22} - M_{12}N_{11} - M_{22}N_{12}}
\end{split}
\end{eqnarray}

\end{appendices}


\begin{thebibliography}{99}

\bibitem{Jirasek-1} M. Jirasek. %
Comparative study on finite elements with embedded discontinuities. %
Comp. Meth. in Appl. Mech. and Eng. 2000; 188: 307--330.

\bibitem{Belytschko-2} T. Belytschko, T. Black. %
Elastic crack growth in finite elements with minimal remeshing. %
Int. J. for Num. Meth. in Eng. 1999; 45(5): 601--620.

\bibitem{Moes-1} N. Moes, J. Dolbow, T. Belytschko. %
A finite element method for crack growth without remeshing. %
Int. J. for Num. Meth. in Eng. 1999; 46(1): 133--150.

\bibitem{Dolbow-1} J. Dolbow, N. Moes, T. Belytschko. %
Discontinuous enrichment in finite elements with a partition of unity method. %
Finite Elements in Analysis and Design 2000; 36(3): 235--260.

\bibitem{Belytschko-3} T. Belytschko, N. Moes, S. Usui, C. Parimi. %
Arbitrary discontinuities in finite elements. %
Int. J. for Num. Meth. in Eng. 2001; 50(4): 993--1013.

\bibitem{Belytschko-4} T. Belytschko, H. Chen, J. X. Xu, G. Zi. %
Dynamic crack propagation based on loss of hyperbolicity and a new discontinuous enrichment. %
Int. J. for Num. Meth. in Eng. 2003; 58: 1873--1905.

\bibitem{Zi-Belytschko} G. Zi, T. Belytschko %
New crack-tip elements for XFEM and applications to cohesive cracks. %
Int. J. Numer. Meth. Engng 2003; 57:2221--2240.

\bibitem{Shiva-1} S. Esna Ashari, S. Mohammadi. %
Delamination analysis of composites by new orthotropic bimaterial extended finite element method. %
Int. J. for Num. Meth. in Eng. 2011; 86: 1507–-1543.

\bibitem{Shiva-2} S. Esna Ashari, S. Mohammadi. %
Fracture analysis of FRP-reinforced beams by orthotropic XFEM. %
Journal of Composite Materials. 2011; 0(0): 1–-23.

\bibitem{Camacho-1} G. T. Camacho, M. Ortiz. %
Computational modeling of impact damage in brittle materials. %
Int. J. of Solids and Structures 1996; 33: 2899--2938.

\bibitem{Ortiz-1} M. Ortiz, A. Pandolfi. %
Finite-deformation irreversible cohesive elements for three-dimensional crack-propagation analysis. %
Int. J. Numer. Meth. Engng. 1999; 44: 1267--1282.

\bibitem{Pandolfi-1} A. Pandolfi, M. Ortiz. %
An efficient adaptive procedure for three-dimensional fragmentation simulations. %
Engineering with Computers 2002; 18(2):148--159.

\bibitem{Cundall-1} P. A. Cundall, O. D. L. Strack. %
A discrete numerical model for granular assemblies. %
Geotechnique 1979; 29: 47--65.

 \bibitem{Hrennikoff-1}  A. Hrennikoff. %
Solution of problems of elasticity by the framework method. %
J Appl Mech 1941; 12: 169-75.

\bibitem{Schlangen-1} E. Schlangen, J.G.M. van Mier. %
Experimental and numerical analysis of micromechanisms of fracture of cement-based composites. %
Cement Concrete Composite 1992; 14:105-118.

\bibitem{Cusatis-1} G. Cusatis, Z. P. \B, L. Cedolin. %
Confinement-shear lattice model for concrete damage in tension and compression: I. Theory. %
J. Eng. Mech. (ASCE) 2003; 129(12): 1439--1448.

\bibitem{Cusatis-2} G. Cusatis, Z. P. \B, L. Cedolin. %
Confinement-shear lattice model for concrete damage in tension and compression: II. Computation and validation. %
J. Eng. Mech. (ASCE) 2003; 129(12): 1449--1458.

\bibitem{Lilliu-1} G. Lilliu, J. G. M. van Mier. %
3D lattice type fracture model for concrete. %
Eng. Fract. Mech. 2003; 70: 927--941.

\bibitem{Cusatis-3} G. Cusatis, Z. P. \B, L. Cedolin. %
Confinement-shear lattice model for fracture propagation in concrete. %
Comput. Methods Appl. Mech. Eng. 2006; 195: 7154--7171.

\bibitem{Bolander-2} J. E. Bolander, S. Saito. %
Fracture analysis using spring network with random geometry. %
Eng. Fract. Mech. 1998; 61(5-6): 569--591.

\bibitem{Bolander-1} J. E. Bolander, K. Yoshitake, J. Thomure. %
Stress analysis using elastically uniform rigid-body-spring networks. %
J. Struct. Mech. Earthquake Eng. (JSCE) 1999; 633(I-49): 25--32.

\bibitem{Bolander-3} J. E. Bolander, G. S. Hong, K. Yoshitake. %
Structural concrete analysis using rigid-body-spring networks. %
J. Comp. Aided Civil and Infrastructure Eng. 2000; 15: 120--133.

\bibitem{Bolander-4} Yip, Mien, Jon Mohle, and J. E. Bolander. %
Automated modeling of three-dimensional structural components using irregular lattices. %
Computer-Aided Civil and Infrastructure Engineering; 2005 20(6): 393--407.

\bibitem{cusatis-ldpm-1} G. Cusatis, D. Pelessone, A. Mencarelli. %
Lattice Discrete Particle Model (LDPM) for failure behavior of concrete. I: Theory. %
Cement and Concrete Composites 2011; 33: 881-890.

\bibitem{cusatis-ldpm-2} G. Cusatis, A. Mencarelli, D. Pelessone, J. Baylot. %
Lattice Discrete Particle Model (LDPM) for failure behavior of concrete. II: Calibration and validation. %
Cement and Concrete Composites 2011; 33: 891-905.

\bibitem{Rezakhani-JMPS}
R. Rezakhani, G. Cusatis, Asymptotic expansion homogenization of discrete fine-scale models with rotational degrees of freedom for the simulation of quasi-brittle materials. \emph{J. Mech. Phys. Solids} 2016; 88: 320--345. 

\bibitem{Leite-1} J.P.B. Leite, V. Slowik, H. Mihashi. %
Computer simulation of fracture processes of concrete using mesolevel models of lattice structures. %
Cement and Concrete Research 2004; 34: 1025-1033.

\bibitem{Donze-1} F. Camborde, C. Mariotti, F.V. Donze. %
Numerical study of rock and concrete behaviour by discrete element modelling. %
Computers and Geotechnics, 2000; 27: 225--247.

\bibitem{Grassl-1} P. Grassl, Z. Bazant, G. Cusatis. %
Lattice-cell approach to quasi-brittle fracture modeling. %
Computational Modelling of Concrete Structures 2006; 930: 263--268 930.

\bibitem{tonti-1} E. Tonti, F. Zarantonello. %
Algebraic formulation of elastostatics: the Cel method. %
Computer Modeling in Engineering and Science 2009; 39(3): 201--236.

\bibitem{Guzey-1} S. G\"{u}zey, B. Cockburn, H. K. and Stolarski. %
The embedded discontinuous galerkin method: application to linear shell problems. %
Int. J. Numer. Meth. Engng. 2007; 70: 757--790.

\bibitem{Shen-1} Y. Shen, A. Lew. %
An optimally convergent discontinuous Galerkin-based extended finite element method for fracture mechanics. %
Int. J. Numer. Meth. Engng. 2010; 82: 716--755.

\bibitem{Abedi-1} R. Abedi, M. A. Hawker, R. B. Haber, K. Matous. %
An adaptive spacetime discontinuous Galerkin method for cohesive models of elastodynamic fracture. %
Int. J. Numer. Meth. Engng. 2010; 81: 1207--1241.

\bibitem{Allman-1} D. J. Allman. %
A compatible triangular element including vertex rotations for plane elasticity analysis. %
Computers and Structures 1984; 19(2): 1--8.

\bibitem{Bergan-1} P. G. Bergan, C. A. Felippa. %
A triangular membrane element with rotational degrees of freedom. %
Comp. Meth. Appl. Mech. Eng. 1985; 50: 25--69.

\bibitem{Zhou-1} X. Zhou, G. Cusatis. %
Tetrahedral finite element with rotational degrees of freedom for Cosserat and Cauchy continuum problems. %
J. Eng. Mech. (ASCE) 2015; 141(2), 06014017.

\bibitem{tonti-2} E. Tonti. %
A direct discrete formulation of field laws: The cell method. %
Computer Modeling in Engineering and Sciences, 2001; 2(2): 237--258.

\bibitem{Polymesh-paulino} C. Talischi, G. H. Paulino, A. Pereira, I. F. Menezes. %
PolyMesher: a general-purpose mesh generator for polygonal elements written in Matlab. %
Structural and Multidisciplinary Optimization, 2012; 45(3): 309--328.

\bibitem{Bazant-CrackBandModel} Z. P. \B,  H. Oh.Byung. %
Crack band theory for fracture of concrete. %
Materiaux et construction, 1983; 16(3): 155-177.


\bibitem{Kalthoff-1} J. F. Kalthoff,  S. Winkler. %
Failure mode transition at high rates of shear loading. %
Int. Conf. on Impact Loading and Dynamic Behavior of Materials 1987; 1: 85--195.

\bibitem{Belytschko-1} T. Belytschko, H. Chen, J. Xu, G. Zi. %
Dynamic crack propagation based on loss of hyperbolicity and a new discontinuous enrichment. %
Int. J. Numer. Meth. Engng 2003; 58: 1873-–1905. 

\bibitem{Song-1} J.-H. Song, H. Wang, T. Belytschko. %
A comparative study on finite element methods for dynamic fracture. %
Comput. Mech. 2008; 42: 239--250.

\bibitem{Xu-1} X.-P. Xu, A. Needleman. %
Numerical simulation of fast crack growth in brittle solids. %
J. Mech. Phys. Solids 1994; 42(9): 1397--1434.

\bibitem{Ramulu-1} M. Ramulu, A. S. Kobayashi. %
Mechanics of crack curving and branching - a dynamic fracture analysis. %
Int. J. Fract. 1985: 27; 187--201.

\bibitem{Ravi-Chandar-1} K. Ravi-Chandar. %
Dynamic fracture of nominally brittle materials. %
Int. J. Fract. 1998: 90; 83--102.


\end{thebibliography}
\end{document}